\documentclass[11pt,twoside]{article}

\font\BBD = msbm10 at 14pt
 \font\bbd = msbm10 at 10pt
\font\bbds = msbm8 at 8pt
\font\bbds = msbm7 at 8pt
\font\csc = cmcsc10 at 10pt
\font\frak = eufm10 at 10pt

\font\xmplbx = cmbx10 scaled \magstep1

 \outer\def\proclaim #1 #2\par{\medbreak
      \noindent{\bf#1\enspace}{\sl#2}\par
      \ifdim\lastskip<\medskipamount \removelastskip\penalty55\medskip\fi}
\def\heading#1#2{%
        \vskip 0pt plus.3\vsize\penalty -10
                          \vskip 0pt plus-.3\vsize\bigskip\medskip\vskip\parskip{\noindent\xmplbx #1.\enspace\ignorespaces#2.\enspace}\ignorespaces\medskip}
\def\proof {\noindent {\it Proof.} \ \ }
\outer\def\remark #1 #2\par{\medbreak
      \noindent{\it#1\enspace}{\rm#2}\par
      \ifdim\lastskip<\medskipamount \removelastskip\penalty55\medskip\fi}

\def\sqr#1#2{{\vcenter{\hrule height.#2pt
              \hbox{\vrule width.#2pt height#1pt \kern#1pt \vrule width.#2pt}
                    \hrule height.#2pt}}}
\def\qe{\mathchoice\sqr54\sqr54\sqr{2.1}3\sqr{1.5}3}

\def\qed{\quad\qe}

\def\C{\mathord{\hbox{\bbd C}}}
\def\c{\mathord{\hbox{\bbds C}}}

\def\A{\mathord{\hbox{\frak a}}}

\def\G{\mathord{\hbox{\frak g}}}

\def\Z{\mathord{\hbox{\bbd Z}}}
\def\ZZ{\mathord{\hbox{\BBD Z}}}

\def\setbrackets #1 #2{\bigr\{\,#1\bigm|#2\,\bigr\}}

\def\dd{\mathord{\hbox{$\partial$}}}

\begin{document}
\thispagestyle{empty}
\thispagestyle{empty}

\leftline{\copyright~ 1998 International Press}
\leftline{Adv. Theor. Math. Phys. {\bf 2} (1998) 1141-1182}  
\vspace{0.4in}

\centerline {\huge \bf Generalized Spencer Cohomology} 
\medskip
\centerline{\huge \bf  and Filtered Deformations of $\ZZ$-graded}
\medskip
\centerline {\huge \bf  Lie Superalgebras}
\vskip 0.2in 
\renewcommand{\thefootnote}{\fnsymbol{footnote}}
\centerline {\bf{$^{a}$Shun-Jen Cheng}\footnote{ partially supported  
by NSC grant 87-2115-M006-005 of the ROC}
 and \bf{$^{b}$Victor G. Kac}\footnote{partially supported  
by NSF grant DMS-9622870}}
\vskip 0.2in 
\centerline{$^{a}$Department of Mathematics,} 
\centerline{National Cheng-Kung University}
\centerline{Tainan,Taiwan}
\centerline{\tt chengsj@mail.ncku.edu.tw}
%\hfill\break\indent\phantom{$^\dagger$}chengsj@mail.ncku.edu.tw}
\medskip
\centerline{$^{b}$Department of Mathematics,}
\centerline{MIT,}
\centerline{Cambridge, MA 02139}
\centerline{\tt kac@math.mit.edu}
\renewcommand{\thefootnote}{}
\footnotetext{\small e-print archive: {\texttt http://xxx.lanl.gov/abs/math.RT/9805039}}
\renewcommand{\thefootnote}{\arabic{footnote}}
\vskip 0.2in
\centerline{\bf Abstract}
\medskip
\begin{center}
\parbox[c]{4.5 in}{\small {\hspace{0.2 in}}In this paper we introduce generalized Spencer cohomology for finite depth $\Z$-graded Lie (super)algebras.  We develop a method of finding filtered deformations of such $\Z$-graded Lie (super)algebras based on t
his cohomology.   As an application we determine all simple filtered deformations of certain $\Z$-graded Lie superalgebras classified in [K3], thus completing the last step in the classification of simple infinite-dimensional linearly compact Lie superalg
ebras.}
\end{center}
\newpage
\pagenumbering{arabic}
\setcounter{page}{1142}

\pagestyle{myheadings}
\markboth{\it GENERALIZED SPENCER COHOMOLOGY ...}{\it S. CHENG, V. KAC}
\heading{0}{Introduction}

Spencer cohomology of a $\Z$-graded Lie algebra $\G=\oplus_{j\ge
-1}\G_j$ of depth $1$ is an important tool for the study of
deformations of geometric structures on a manifold [S], [GS], [SS],
[KN].  It works best, however, when the corresponding transitive
pseudogroup of transformations is irreducible, i.e.~admits no
invariant differential systems, integrable or not.  In the primitive
(but reducible) case, i.e.~when there are no integrable
differential systems, it is natural to pick a minimal invariant
(non-integrable) differential system, which leads to $\Z$-graded Lie
algebras of depth $h\ge 1$ [W]: $\G=\oplus_{j\ge -h}\G_j$. This brings
us to the generalized Spencer cohomology (Section 1).

The present paper is a part of the program of classification, up to
formal equivalence, of simple infinite-dimensional Lie superalgebras
of vector fields on a finite-dimensional supermanifold [K3].  Like in
the Lie algebra case, an important ingredient in this classification
is the description of all simple filtered deformations of a given
graded Lie superalgebras cf.~[SS], [KN], [W], [G].

Let $L$ be a linearly compact Lie (super)algebra, that is a complete topological Lie (super)algebra, which admits a fundamental system of neighborhoods of $0$ consisting of subspaces of finite codimension.  (The formal completion of a Lie (super)algebra o
f vector fields on a finite-dimensional (super)manifold $X$ at a neighborhood of a point of $X$ is of this kind.) Provided that $L$ is simple (i.e. has no non-trivial closed ideals), one can construct a filtration of $L$ by open (and hence closed) subspac
es
$$L=L_{-h}\supset L_{-h+1}\supset\cdots\supset L_0\supset L_1\supset\cdots,$$
such that the associated graded Lie (super)algebra ${\rm Gr}L=\oplus_{j=-h}^{\infty}\G_j$, $\G_j=L_j/L_{j+1}$, of depth $h$ has the properties [W]:\newline
 (G0) ${\rm dim}\G_j<\infty$,\newline
 (G1) $\G_{-j}=\G_{-1}^j$, for $j\ge 1$,\newline
 (G2) if $a\in\G_j$, $j\ge 0$, then $[a,\G_{-1}]=0$ implies that $a=0$,\newline (G3) the representation of $\G_0$ on $\G_{-1}$ is irreducible.

In the Lie algebra case such a filtration is unique, provided that ${\rm dim}L=\infty$ [G], and it is not too hard to classify all $\Z$-graded Lie algebras satisfying properties (G0)--(G3) (see [K1] or [G]).  However, in the Lie superalgebra case there ar
e many such filtrations and it is all but impossible to classify all $\Z$-graded Lie superalgebras satisfying (G0)--(G3).  The basic idea of [K3] is to choose a ``maximally even'' $L_0$; then the representation of $\G_0$ on $\G_{-1}$ satisfies much more s
evere restrictions than (G3) (cf.~[G]), which makes it possible to classify such $\Z$-graded Lie superalgebras.

The next step is to describe, for each $\Z$-graded Lie (super)algebra $\G$ of the obtained list, all simple {\it filtered deformations} of $\G$, i.e.~all simple linearly compact Lie (super)algebras $L$ such that ${\rm Gr}L\cong\G$.  Of course, if $\G$ is 
a simple $\Z$-graded Lie (super)algebra, then its completion $\bar{\G}$ in topology defined by the fundamental system $\G_{(k)}=\oplus_{i\ge k}\G_i$, $k\in\Z_+$, is a simple filtered deformation, called the {\it trivial filtered deformation}.  It is easy 
to show (cf.~Corollary 2.2) that if $\G_0$ contains a non-zero central element, then $\G$ has only a trivial deformation.

In the Lie algebra case the only remaining examples are the two series of $\Z$-graded Lie algebras of depth $1$, which consist of divergence free and Hamiltonian vector fields with polynomial coefficients.  In these two cases one can either use the classi
cal Spencer cohomology as in [SS], [KN], or some more ``pedestrian'' arguments, as in [W], [K2], [K3], to show that all filtered deformations are trivial.

However, in the Lie superalgebra case there are many more cases of $\Z$-graded Lie superalgebras, and only for some of them the ``pedestrian'' arguments work (cf.~[K3]).  Also, we do not have at our disposal a Serre type vanishing theorem for Spencer coho
mology as in the Lie algebra case (cf.~[KN]).  Moreover, there are several series of $\Z$-graded Lie superalgebras of depth $h\ge 2$ to which the classical Spencer cohomology is not applicable.

The aim of the present paper is to show how to resolve these difficulties.  In Section 1 we introduce {\it generalized} Spencer cohomology, which is applicable to graded Lie superalgebras of arbitrary depth $h$.  In Section 2 we show that filtered deforma
tions are described by the invariant Spencer $2$-cocycles, provided that $\G$ is an {\it almost full prolongation}.  We introduce the latter notion since, unlike in the Lie algebra case, not all $\Z$-graded algebras in question are full prolongations (mea
ning that the first Spencer cohomology is trivial), but all, except for one of them, happen to be almost full prolongations.

After describing in Section 3 all examples of $\Z$-graded Lie superalgebras determination of whose filtered deformations was left out in [K3], we apply to them in Section 4 the techniques developed in Section 2.  We find that, unlike in the Lie algebra ca
se, there are three series of $\Z$-graded Lie superalgebras that do admit a (unique) non-trivial filtered deformation (Theorems 5.1 and 5.2), and the rest do not (Theorems 4.1--4.4).  Note that one of these filtered deformations was discovered by Kotchetk
off [Ko], and that the two filtered deformations discovered in this paper are isomorphic.

We would like to thank Yuri Kotchetkoff for very useful correspondence.

All vector spaces, algebras and tensor products in this paper are considered over the field $\C$ of complex numbers.

\heading{1}{Generalized Spencer cohomology}

Let $\G=\oplus_{j=-h}^{\infty}\G_j$ be a Lie superalgebra with a $\Z$-gradation (compatible with its $\Z_2$-gradation) of finite depth $h$, where $h$ is a positive integer.  We have $[\G_i,\G_j]\subset\G_{i+j}$, and we shall always assume that ${\rm dim}\
G_j<\infty$ for all $j\ge-h$.

Set $\G_-=\oplus_{j=-h}^{-1}\G_j$.  Obviously $\G_-$ is a (finite-dimensional) subalgebra of $\G$ and hence $\G_-$ acts on $\G$ via the adjoint representation, so that we may consider $H^*(\G_-;\G)$, the cohomology groups of $\G_-$ with coefficients in it
s adjoint representation in $\G$.  Recall (see e.g.~[F]) that the space of cochains is
$$C^*(\G_-;\G)=\Lambda^{*}(\G_-)\otimes\G,$$
where the exterior product $\Lambda^*$ is understood in the usual super sense. The space of the $j$-cochains ($j\in\Z_+$) is then
$$C^{j}(\G_-;\G)=\oplus_{-1\ge i_1\ge i_2\ge\cdots\ge i_j\ge -h}(\G_{i_1}^*\wedge\cdots\wedge\G_{i_j}^*)\otimes\G,$$
on which the coboundary operator $d$ acts as
($p+q=j+1$)

\begin{eqnarray}
&&(dc)(x_1,\cdots,x_p,y_1,\cdots,y_q) \nonumber= \\&& \sum_{1\le s<t\le p}(-1)^{s+t-1} c([x_s,x_{t}],x_1,\cdots,\hat{x}_s,\cdots,\hat{x}_{t},\cdots,x_p,y_1,\cdots,y_q)\nonumber \\
&&+ \sum_{s=1}^{p}\sum_{t=1}^q (-1)^{s-1}c(x_1,\cdots,\hat{x}_s,\cdots,x_p,[x_s,y_{t}],y_1,\cdots,\hat{y}_{t},\cdots,y_q))\nonumber \\
&&+ \sum_{1\le s<t\le q} c([y_s,y_{t}],x_1,\cdots,x_p,y_1,\cdots,\hat{y}_s,\cdots,\hat{y}_{t},\cdots,y_q)\nonumber \\
&&+ \sum_{s=1}^p (-1)^s [x_s,c(x_1,\cdots,\hat{x}_s,\cdots,x_p,y_1,\cdots,y_q)]\nonumber \\
&&+ (-1)^{p-1}\sum_{s=1}^q [y_s,c(x_1,\cdots,x_p,y_1,\cdots,\hat{y}_s,\cdots,y_q)],\nonumber 
\end{eqnarray}
where $x_1,\cdots,x_p\in {(\G_-)}_{\bar{0}}$ and $y_1,\cdots,y_q\in {(\G_-)}_{\bar{1}}$. We have then $H^*(\G_-;\G)={\rm Ker}d/{\rm Im}d$.

Note that $C^j(\G_-;\G)$ is $\Z$-graded by letting ${\rm deg}\G_i=-{\rm deg}\G^*_i=i$.  This gradation induces a $\Z$-gradation on $H^j(\G_-;\G)$:
$$H^j(\G_-;\G)=\oplus_{l\in\Z}H^{l,j}(\G_-;\G),$$
where $H^{l,j}(\G_-;\G)$ denotes the $l$-th graded component of $H^j(\G_-;\G)$.

We will call the vector space $H^{l,j}(\G_-;\G)$ the {\it $(l,j)$-th (generalized) Spencer cohomology group} of the $\Z$-graded Lie superalgebra $\G$ and we will call elements in  ${\rm Ker}d$ (respectively ${\rm Im}d$) {\it Spencer cocycles} (respectivel
y {\it Spencer coboundaries}).  In this paper only those $H^{l,j}(\G_-;\G)$, for which $l\ge 0$ will play a role.  This definition is a generalization of the classical Spencer cohomology defined for $h=1$.  We would like to point out that the classical Sp
encer cohomology $H^{p,q}$ (cf.~[Sp]) would in our definition correspond to $H^{p+q-1,q}$.

From the definition it is obvious that $H^{*,0}$ is the subspace of \newline
$\G_-$-invariants in $\G$ so that we have
$$H^{k,0}(\G_-;\G)=(\G_k)^{\G_-}.$$
The Lie superalgebra $\G$ is called {\it transitive} if $H^{k,0}(\G_-;\G)=0$ for all $k\ge 0$, i.e.~$\G$ is transitive if the conditions $[\G_-,a]=0$ and $a\in\oplus_{j\ge 0}\G_j$ imply that $a=0$.

\remark{Remark 1.1.} The transitivity property is equivalent to (G2), provided that (G1) holds (see Introduction).

A linear map $\alpha:\G_-\rightarrow\G$ is called a {\it derivation} of $\G_-$ into $\G$ if for all $x,y\in\G_-$ we have $\alpha([x,y])=[\alpha(x),y]+(-1)^{p(\alpha)p(x)}[x,\alpha(y)]$. Evidently, the space of all derivations ${\rm der}_{\c}(\G_-,\G)$ is 
$\Z$-graded so that we may write ${\rm der}_{\c}(\G_-,\G)=\oplus_{l\in\Z}{\rm der}_{\c}(\G_-,\G)_l$.  Furthermore every element of $\G$ itself defines a derivation of $\G_-$ into $\G$.  It follows from the definition that $H^{*,1}(\G_-;\G)={\rm der}_{\c}(
\G_-,\G)/\G$ so that
$$H^{l,1}(\G_-;\G)={\rm der}_{\c}(\G_-,\G)_l/\G_l.$$

Let $\G_{\le 0}$ denote the subalgebra $\G_-\oplus\G_0$ of $\G$.  We will say that $\G$ is a {\it full prolongation of $\G_{\le 0}$ of degree k}, if $\G$ contains all derivations of $\G_-$ into $\G$ of degree $\ge k$.  This is equivalent to saying that $H
^{l,1}(\G_-;\G)=0$, for $l\ge k$.  Note that a full prolongation of $\G_{\le 0}$ of degree $1$ is uniquely determined (since in this case $\G_j$ is just ${\rm der}_{\c}(\G_-,\G)_j$, for $j\ge 1$); in this case we shall call $\G$ the {\it full prolongation
} of $\G_{\le 0}$.

\remark{Remark 1.2.} For most cases the notion of full prolongation is adequate for the study of filtered deformations.  However, in some cases, full prolongation is too strong an assumption, and should be replaced by a weaker notion, which we shall call 
an {\it almost full prolongation}.  We shall take this up in the next section after introducing filtered deformations, which turn out to be closely related to $H^{*,2}(\G_-;\G)$.
\heading{2}{Preliminaries on filtered deformations}

Let $L$ be a filtered Lie superalgebra of finite depth $h$, where $h$ is a positive integer.  This means that $L$ is a Lie superalgebra with a sequence of subspaces (compatible with the $\Z_2$-gradation of $L$)
$$L=L_{-h}\supset L_{-h+1}\supset \cdots L_{-1}\supset L_0\supset L_{1}\cdots\supset L_n\supset \cdots,$$
such that $[L_i,L_j]\subset L_{i+j}$.  We shall assume in this paper that ${\rm dim}_{\c}L_{j}/L_{j+1}<\infty$, for all $j$.  The filtration in a natural way induces a topology on $L$.  The condition $[L_i,L_j]\subset L_{i+j}$ makes $L$ into a topological
 Lie superalgebra.  We will say that $L$ is {\it complete}, if $L$ is complete with respect to this topology.  In this paper we shall always deal with complete filtered Lie superalgebras. Let $\G=\oplus_{j=-h}^{\infty}\G_j$, where $\G_j=L_j/L_{j+1}$, be i
ts associated graded Lie superalgebra.  We let $\G_{(j)}=\oplus_{i\ge j}\G_i$.  This defines a filtration on $\G$.  The completion of $\G$ with respect to the topology induced by this filtration will be denoted by $\bar{\G}$.

For each $j\ge-h$ we may choose a subspace $V_j$ of $L_j$ so that $V_j\oplus L_{j+1}=L_j$ as vector spaces.  We may identify $V_j$ with $\G_j$  so that in the vector space $\bar{\G}=\prod_{j}\G_j=\prod_{j}V_j=L$ we may define two Lie brackets.  Namely, $[
\cdot,\cdot]$, which is the Lie bracket of the Lie superalgebra $\bar{\G}$, and $[\cdot,\cdot]_{1}$, which is the Lie bracket of the Lie superalgebra $L$.  We have for $x,y\in\G=\oplus_{j}\G_j$:
$$[x,y]_{1}=[x,y]+\sum_{i\ge 1}\mu_{i}(x,y),\eqno{(2.1)}$$
where $\mu_i:\G\wedge\G\rightarrow\G$ is an even super-skewsymmetric bilinear map such that $\mu_i(\G_j\wedge\G_s)\subset\G_{j+s+i}$ for each $i=1,2,\cdots$.  Note that for each $\epsilon\in\C^*$ the map $\varphi_{\epsilon}:\bar{\G}\rightarrow\bar{\G}$, d
efined by $\varphi_{\epsilon}(x)=\epsilon^jx$, if $x\in\G_j$, is a continuous automorphism of the Lie superalgebra $\bar{\G}$, provided that $\epsilon\not=0$.  Applying $\varphi_{\epsilon}$, with $\epsilon\not=0$, to both sides of (2.1) and dividing by an
 appropriate power of $\epsilon$, we obtain, letting $[x,y]_{\epsilon}=\varphi_{\epsilon}([x,y]_1)$:
$$[x,y]_{\epsilon}=[x,y]+\sum_{i\ge 1}\mu_{i}(x,y){\epsilon}^i.\eqno{(2.2)}$$
The bracket $[x,y]_{\epsilon}$ defines a Lie superalgebra structure on the space $\bar{\G}$.  If $\epsilon\not=0$, the obtained Lie superalgebra, which we denote by $\bar{\G}_{\epsilon}$, is isomorphic to $\bar{\G}_1$.  If $\epsilon=0$, it is isomorphic t
o $\bar{\G}$.
We will sometimes call $[\cdot,\cdot]_{\epsilon}$ the {\it deformed bracket} of $[\cdot,\cdot]$, and $\bar{\G}_{\epsilon}$ with a deformed bracket (or $L$) a {\it filtered deformation} of $\G$. A filtered deformation is said to be a {\it trivial} deformat
ion, if it is isomorphic to $\bar{\G}$.

We have associated to a filtered deformation $\bar{\G}_{\epsilon}\cong L$ of $\G$ a sequence of bilinear maps $\mu_i:\G\wedge\G\rightarrow\G$, $i=1,2,\cdots$.  We shall call the sequence $\{\mu_1,\mu_2,\cdots\}$ a {\it defining sequence of this filtered d
eformation}.

\remark{Remark 2.1.} Of course, a different choice of the subspaces $V_j$ gives rise to a different defining sequence in general.  Hence a filtered deformation may be represented by different defining sequences.    Thus we may study a filtered deformation
 by analyzing the effect of a different choice of the subspaces $V_j$  on the resulting defining sequence.  Clearly, a filtered deformation is trivial if and only if we may choose the subspaces $V_j$ in such a way that the resulting defining sequence cons
ists of zero maps.

Let $x$, $y$ and $z$ be homogeneous (both in the $\Z$- and $\Z_2$-grading) elements of $\G$. The Jacobi identity in $\bar{\G}_{\epsilon}$ gives
$$[x,[y,z]_{\epsilon}]_{\epsilon}=[[x,y]_{\epsilon},z]_{\epsilon}+p(x,y)[y,[x,z]_{\epsilon}]_{\epsilon},\quad {\rm where\ }p(x,y)=(-1)^{p(x)p(y)}.$$
Substituting (2.2) into this expression gives an identity in power series in ${\epsilon}$ with coefficients in $\G$.  We collect the coefficient of ${\epsilon}^k$ and obtain the following identity in $\G$:
\begin{eqnarray}
&&[x,\mu_k(y,z)]+\mu_k(x,[y,z])+\sum_{1\le i,j<k}\mu_i(x,\mu_j(y,z))\nonumber\\
&=&\mu_k([x,y],z) +[\mu_k(x,y),z]
+ \sum_{1\le i,j<k}\mu_i(\mu_j(x,y),z) +p(x,y)[y,\mu_k(x,z)]\nonumber\\  &+& p(x,y)\mu_k(y,[x,z])
+ p(x,y)\sum_{1\le i,j<k}\mu_i(y,\mu_j(x,z)).\;\;\;\;\;\;\;\;\;\;\;\;\;\;\;\;\;\;\;\;\;\;\;\;\;{(2.3)}\nonumber
\end{eqnarray}

\proclaim{Proposition 2.1.} The first non-zero term $\mu_k$ in (2.2) is an even $2$-cocycle of $\G$ with coefficients in the adjoint representation.

\proof Since $\mu_i=0$ for all $i<k$, from (2.3) we get
\begin{eqnarray}
[x,\mu_k(y,z)]+\mu_k(x,[y,z])&=&\mu_k([x,y],z)+[\mu_k(x,y),z]\nonumber \\
&+& p(x,y)[y,\mu_k(x,z)] \nonumber \\
&+& p(x,y)\mu_k(y,[x,z]).\quad\quad\quad\quad\quad\quad\quad\quad\quad{(2.4)}\nonumber
\end{eqnarray}
But this precisely means that $\mu_k$ is a $2$-cocycle.  {$\qed$\medskip}

We rewrite (2.4) as
$$\mu_k(x,[y,z])-[\mu_k(x,y),z]-p(x,y)[y,\mu_k(x,z)]=$$
$$\mu_k([x,y],z)+p(x,y)\mu_k(y,[x,z])-[x,\mu_k(y,z)].$$
The right hand side above is precisely
$$x\cdot\mu_k(y,z),$$
while the left hand side is $-df^k_x(y,z)$, where $f^k_x:\G\rightarrow\G$ is given by $f^k_x(y)=\mu_k(x,y)$ and $d$ is the coboundary operator.  Thus it follows from the proof of Proposition 2.1 that (2.3) may be rewritten as
\begin{eqnarray}
x\cdot\mu_k(y,z)+df_x(y,z)&=&-\sum_{1\le i,j<k}\mu_i(x,\mu_j(y,z))
+\sum_{1\le i,j<k}\mu_i(\mu_j(x,y),z)\nonumber \\
&+& p(x,y)\sum_{1\le i,j<k}\mu_i(y,\mu_j(x,z)),\quad\quad\quad\quad\quad\quad (2.5)\nonumber
\end{eqnarray}
 Of course (2.3) can also be rewritten as
\begin{eqnarray}
d\mu_k(x,y,z)&=&-\sum_{1\le i,j<k}\mu_i(x,\mu_j(y,z))
+\sum_{1\le i,j<k}\mu_i(\mu_j(x,y),z)\nonumber \\
&+&p(x,y)\sum_{1\le i,j<k}\mu_i(y,\mu_j(x,z)).\quad\quad\quad\quad\quad\quad\quad\quad\quad\quad\quad{(2.6)}\nonumber 
\end{eqnarray}

Here is the key observation (due to Kobayashi and Nagano [KN] in the case $h=1$):

\proclaim{Proposition 2.2.} The first non-zero term $\mu_k$ in (2.2) restricted to $\G_-$ defines an even $\G_0$-invariant element in $H^{k,2}(\G_-;\G)$.

\proof By Proposition 2.1 $\mu_k|_{\G_-\times\G_-}$ is a $2$-cocycle.  Now (2.5) with $x\in\G_0$ and $y,z\in\G_-$ means precisely that it is $\G_0$-invariant in $H^*(\G_-;\G)$ (since its right-hand side is zero).  {$\qed$\medskip}

\remark{Remark 2.2.} Proposition 2.2 also follows from Proposition 2.1 as follows: $\mu_k$ defines a $\G$-invariant element of $H^*(\G_-;\G)$, hence it is $\G_0$-invariant in $H^*(\G_-;\G)$, when restricted to $\G_-$.

\proclaim{Proposition 2.3.} Let $\bar{\G}_{\epsilon}$ and $\bar{\G}_{\epsilon}'$ be two filtered deformations given by defining sequences $\{\mu_1,\mu_{2},\cdots\}$ and $\{\mu'_1,\mu'_{2},\cdots\}$, respectively.  Suppose that $(\mu_k-\mu'_k)|_{\G_-\times
\G_-}$ is a Spencer coboundary for some $k\ge 1$.  Then $\bar{\G}'_{\epsilon}$ has a defining sequence $\{\mu_1'',\mu_2'',\cdots\}$ such that $\mu''_i=\mu'_i$, for $i<k$, and $\mu_k''|_{\G_-\times\G_-}=\mu_k|_{\G_-\times\G_-}$.  (In other words, one can c
hange the defining sequence of $\bar{\G}'_{\epsilon}$ such that its first $k-1$ terms are unchanged, and its $k$-th term becomes the $k$-th term of the defining sequence of $\bar{\G}_{\epsilon}$ when restricted to $\G_-$.)

\proof We have for $x,y\in\G_-$
$$[x,y]_{\epsilon}=[x,y]+\sum_{i\ge 1}\mu_i(x,y){\epsilon}^i,\quad [x,y]'_{\epsilon}=[x,y]+\sum_{i\ge 1}\mu'_i(x,y){\epsilon}^i.$$
By assumption $\mu_k-\mu_k'$ is a Spencer coboundary, hence there exists an $f:\G_-\rightarrow\G$ such that $df=\mu_k-\mu_k'$.  We define an injective map $\rho_k:\G_-\rightarrow\G$ via
$$\rho_k(x)=x+f(x)\epsilon^k,\quad\forall x\in\G_-.$$
A simple calculation shows that
$$[\rho_k(x),\rho_k(y)]'_{\epsilon}=\rho_k([x,y])+\sum_{i<k}\mu_i'(x,y){\epsilon}^i+\mu_k(x,y)\epsilon^k+\sum_{i>k}\mu_i''(x,y){\epsilon}^i.$$
But now obviously $\rho_k(\G_-)+\G_{(0)}=\G$.  Thus replacing $\G_-$ by $\rho_k(\G_-)$ in $\bar{\G}_{\epsilon}'$, which correspond to a new choice of $V_j$, for $j<0$, we obtain a defining sequence with the desired property. {$\qed$\medskip}

Combining Propositions 2.2 and 2.3 we obtain

\proclaim{Corollary 2.1.} Let $\bar{\G}_{\epsilon}$ be a filtered deformation of a graded Lie superalgebra $\G$ with defining sequence $\{\mu_1,\mu_2,\cdots\}$. Suppose that $H^{j,2}(\G_-;\G)$ has no non-trivial even $\G_0$-invariant vectors for any $j\ge
 1$.  Then $\bar{\G}_{\epsilon}$ has a defining sequence $\{\mu_1',\mu_2',\cdots\}$ such that $\mu_{j}'|_{\G_-\times\G_-}$ is identically zero for all $j\ge 1$.

\proof Since $\mu_1|_{\G_-\times\G_-}$ is a Spencer coboundary by Proposition 2.1, the filtered deformation $\bar{\G}_{\epsilon}$ has a defining sequence such that $\{\mu_1',\mu_2',\cdots\}$, where $\mu_1'|_{\G_-\times\G_-}=0$ by Proposition 2.3.  Now  $\
mu_2'|_{\G_-\times\G_-}$ is a Spencer coboundary, hence $\bar{\G}_{\epsilon}$ has a defining sequence $\{\mu_1',\mu_2'',\cdots\}$, such that $\mu_2''|_{\G_-\times\G_-}=0$.  Repeating this procedure, we may make $\mu_3'''|_{\G_-\times\G_-}=0$ as well, etc.
  Since $\bar{\G}_{\epsilon}$ is complete, we may take the limit. {$\qed$\medskip}

Let $\G=\oplus_{j=-h}^{\infty}\G_j$ be a $\Z$-graded Lie superalgebra.  Suppose that $\bar{\G}_{\epsilon}$ is a filtered deformation of $\G$ with defining sequence $\{0,\cdots,0,\mu_k,\mu_{k+1},\cdots\}$.  Let $\A$ be a maximal reductive subalgebra of $\G
_{\bar{0}}$ and suppose that either $\A$ is semisimple or has a $1$-dimensional center $\C c$, where ${\rm ad}c$ acts on $\G_{j}$ as the scalar $j$ for each $j\in\Z$.  By Proposition 2.1, $\mu_k|_{\A\times\A}$ is a $2$-cocycle of $\A$ with coefficients in
 the $\A$-module $\G_k$.  Due to our assumptions on $\A$, by Whitehead's second lemma $\mu_k|_{\A\times\A}$ is a coboundary.  As in the proof of Corollary 2.1 we may find a defining sequence $\{0,\cdots,0,\mu_k',\mu_{k+1}',\cdots\}$ for $\bar{\G}_{\epsilon}$ such that $\mu'_j|_{\A\times\A}$ is identically zero for all $j$.  Thus we may assume that
$$[a,b]_{\epsilon}=[a,b],\quad\forall a,b\in\A.$$

Note that $\mu_k:\G_0\otimes\G_j\rightarrow\G_{j+k}$ induces a map $\nu_{k}|_{\A}:\A\rightarrow\G_j^*\otimes\G_{j+k}$.  It is easy to show, using the Jacobi identity and the fact that $\mu_k|_{\A\times\A}$ is identically zero, that $\nu_k|_{\A}$ is a $1$-
cocycle of $\A$ with coefficients in $\G_j^*\otimes\G_{j+k}$.  If $\A$ is semisimple, then by Whitehead's first lemma $\nu_k|_{\A}$ is a coboundary.  If $\A$ has a non-zero center, then it acts non-trivially on $\G_j^*\otimes\G_{j+k}$ and hence $\nu_k|_{\
A}$ is a coboundary as well.  Arguing as before, we may assume that for all $j$
$$[a,x]_{\epsilon}=[a,x],\quad\forall a\in\A,\forall x\in\G_j.$$

Now let $x\in\G_s$ and $y\in\G_l$.  Taking bracket in $\bar{\G}_{\epsilon}$ of an element $a\in\A$ with $[x,y]_{\epsilon}$ we obtain
\begin{eqnarray}
[a,[x,y]_{\epsilon}]_{\epsilon}&=&[a,[x,y]]_{\epsilon}+\sum_{i=k}^{\infty}[a,\mu_i(x,y)]_{\epsilon}{\epsilon}^i\nonumber \\
&=&[a,[x,y]]+\sum_{i=k}^{\infty}[a,\mu_i(x,y)]{\epsilon}^i.\nonumber 
\end{eqnarray}

On the other hand, by Jacobi identity in $\bar{\G}_{\epsilon}$ the same quantity is equal to
\begin{eqnarray}
&&[[a,x]_{\epsilon},y]_{\epsilon}+[x,[a,y]_{\epsilon}]_{\epsilon}\nonumber \\
&=&[[a,x],y]_{\epsilon}+[x,[a,y]]_{\epsilon}\nonumber \\
&=&[a,x],y]+\sum_{i=k}^{\infty}\mu_i([a,x],y){\epsilon}^i+[x,[a,y]]+\sum_{i=k}^{\infty}\mu_i(x,[a,y]){\epsilon}^i.\nonumber 
\end{eqnarray}
Comparing the coefficients of ${\epsilon}^i$ we obtain
$$[a,\mu_i(x,y)]=\mu_i([a,x],y)+\mu_i(x,[a,y]),$$
which means precisely that the map $\mu_i|_{\G_s\times\G_t}:\G_s\otimes\G_t\rightarrow\G_{s+t+i}$
is a homomorphism of $\A$-modules for every $i\ge k$ and $s,t\ge-h$.  
We thus have proved the following proposition:

\proclaim{Proposition 2.4.} Let $\G=\oplus_{j=-h}^{\infty}\G_j$ be a $\Z$-graded Lie superalgebra and let $\A\subset\G_0$ be a maximal reductive subalgebra of $\G_{\bar{0}}$.  Suppose that either $\A$ is semisimple or the center of $\A$ is $\C c$, where $
{\rm ad}c$ acts on $\G_{j}$ as $j$, for every $j\in\Z$.  Then every filtered deformation of $\G$ has a defining sequence $\{\mu_1,\mu_2,\cdots\}$ such that $\mu_i(\A,\G)=0$ and $\mu_i:\G_s\otimes\G_t\rightarrow\G_{s+t+i}$ is a homomorphism of $\A$-modules
, for $i=1,2,\cdots$.

In other words, Proposition 2.4 says that in every filtered deformation $L$ of $\G$ one can choose a subalgebra $\A'\subset L_0$ which maps isomorphically to $\A$ under the map $L_0\rightarrow\G_0$ and one can choose a subspace $V_j$ in $L_j$ complementar
y to $L_{j+1}$ for each $j\ge -h$ such that $\A'\subset V_0$ and $[\A',V_j]\subset V_j$.  From this we obtain immediately the following (well-known) corollary, which takes care of the case when $\A$ has a non-trivial center.

\proclaim{Corollary 2.2.} Let $\G=\oplus_{j=-h}^{\infty}\G_j$ be a graded Lie superalgebra of depth $h$.  Suppose that $\G_0$ contains an element $c$ such that ${\rm ad}c|_{\G_j}=j$.  Then $\G$ has no non-trivial filtered deformations.

From now on we shall assume that $\G$ is transitive.

\proclaim{Proposition 2.5.} Let $\G=\oplus_{j=-h}^{\infty}\G_j$ be a transitive graded Lie superalgebra.  Suppose that $\{\mu_1,\mu_2,\cdots\}$ is a defining sequence of a filtered deformation $\bar{\G}_{\epsilon}$ of $\G$.  Then $\mu_i$ is completely det
ermined by its restriction to ${\G_-\times\G}$.

\proof Let $a,b\in\G_0$ and $x\in\G_-$.  We have
\begin{eqnarray}
[x,[a,b]_{\epsilon}]_{\epsilon}&=&[x,[a,b]]_{\epsilon}+\sum_{i=1}^{\infty}[x,\mu_i(a,b)]_{\epsilon} {\epsilon}^i\nonumber \\
&=&[x,[a,b]]+[x,\mu_1(a,b)]{\epsilon}+\mu_1(x,[a,b]){\epsilon}+o({\epsilon}^2).\nonumber 
\end{eqnarray}

Now obviously
\begin{eqnarray}
[x,[a,b]_{\epsilon}]_{\epsilon}&=&[[x,a]_{\epsilon},b]_{\epsilon}+p(x,a)[a,[x,b]_{\epsilon}]_{\epsilon}\nonumber \\
&=&[[x,a],b]+p(x,a)[a,[x,b]]+\mu_1([x,a],b){\epsilon}+[\mu_1(x,a),b]{\epsilon}\nonumber \\
&+&p(x,a)[a,\mu_1(x,b)]{\epsilon}+p(x,a)\mu_1(a,[x,b]){\epsilon}+o({\epsilon}^2).\nonumber 
\end{eqnarray}
Hence 
\begin{eqnarray}
[x,\mu_1(a,b)]+\mu_1(x,[a,b])&=&\mu_1([x,a],b)+[\mu_1(x,a),b]\nonumber \\
&+&p(x,a)[a,\mu_1(x,b)]+p(x,a)\mu_1(a,[x,b]).\;\;{(2.7)}\nonumber 
\end{eqnarray}
By assumption $\mu_1|_{\G_-\times\G}$ is known.  Hence the only term in (2.7) that is not determined is $\mu_1(a,b)$.  However, since $\G$ is transitive, $\mu_1(a,b)$ is uniquely determined by (2.7).  Thus $\mu_1|_{\G_0\times\G_0}$ is determined by $\mu_1
|_{\G_-\times\G}$.

Now suppose that $a\in\G_0$ and $b\in\G_k$.  We will argue inductively.  Suppose that $\mu_1|_{\G_-\times\G},\mu_1|_{\G_0\times\G_0},\cdots,\mu_1|_{\G_0\times\G_{k-1}}$ are uniquely determined.  From (2.7) again we see that the only term that is not deter
mined is $\mu_1(a,b)$.  By transitivity again $\mu_1(a,b)$ must be uniquely determined.  Hence $\mu_1|_{\G_-\times\G}$ and $\mu_1|_{\G_0\times\G}$ are uniquely determined.

Now suppose that $a,b\in\G_1$.  Again from (2.7) and transitivity we see that $\mu_1(a,b)$ is uniquely determined.  Similarly $\mu_1|_{\G_1\times\G}$ is uniquely determined.  Proceeding this way we see that $\mu_1$ is uniquely determined by $\mu_1|{\G_-\times\G}$.

Now $\mu_2$ satisfies equation (2.7) up to a function depending only on $\mu_1$ by (2.3).  Since $\mu_1$ is already uniquely determined, we may proceed as before to show that $\mu_2$ is uniquely determined by $\mu_2|_{\G_-\times\G}$ and $\mu_1$.  Similarl
y $\mu_3$ satisfies equation (2.7) up to a function depending on $\mu_1$ and $\mu_2$.  Hence $\mu_3$ is uniquely determined by $\mu_3|_{\G_-\times\G}$, $\mu_1$ and $\mu_2$ etc.  This completes the proof. {$\qed$\medskip}

\proclaim{Proposition 2.6.} Let $\G=\oplus_{j=-h}^{\infty}\G_j$ be a transitive graded Lie superalgebra such that $\{\mu_1,\mu_2,\cdots\}$ and $\{\mu'_1,\mu_2',\cdots\}$ define two filtered deformations of $\G$.  Suppose that $\mu_i|_{\G_-\times\G}=\mu'_i
|_{\G_-\times\G}$ for $i<k$ and $\mu_k|_{\G_-\times\G_-}=\mu'_k|_{\G_-\times\G_-}$ for some $k\ge 1$.  Assume that $\G$ is a full prolongation of $\G_{\le 0}$ of degree $k$. Then $\bar{\G}'_{\epsilon}$ has a defining sequence $\{\mu_1'',\mu_2'',\cdots\}$ 
such that $\mu''_i=\mu_i$, for $i=1,\cdots,k$.  Furthermore $\mu_i''|_{\G_-\times\G_-}{=}\mu'_i|_{\G_-\times\G_-}$ for all $i$. (In other words $\bar{\G}_{\epsilon}$ and $\bar{\G}'_{\epsilon}$ have defining sequences that coincide up to the $k$-th term an
d coincide when restricted to $\G_-\times\G_-$.)

\proof Let $a\in\G_0$ and $x\in\G_-$.
By Proposition 2.5 it follows that $\mu_i=\mu_i'$ for $i<k$.  Now from this, the fact that $(\mu_k-\mu'_k)|_{\G_-\times\G_-}$ is identically zero and (2.5) it is easy to see that the map $f_a^k:\G_-\rightarrow\G$ defined by
$$f_a^k(x):=\mu_k(a,x)-\mu'_k(a,x),\quad x\in\G_-$$
is a Spencer $1$-cocycle.  Hence by hypothesis there exists an element $v_a\in\G$ such that $f_a^k(x)=[v_a,x]$, for all $x\in\G_-$.  Now set $\rho_0^k(a)=a-v_a\epsilon^k$, for all $a\in\G_0$.  It follows that
\begin{eqnarray}
[\rho_0^k(a),x]'_{\epsilon}&=&[a,x]+\sum_{i<k}\mu'_i(a,x){\epsilon}^i\nonumber \\
&+&\mu_k(a,x)\epsilon^k+\sum_{i>k}\mu''_i(a,x){\epsilon}^i,\forall a\in\G_0,\forall x\in\G_-.\nonumber 
\end{eqnarray}

Next let $b\in\G_1$ and $x\in\G_-$.
Using the fact that $(\mu_k-\mu'_k)$, restricted to ${\G_-\times\G_-}$ and $\G_-\times\G_0$, is identically zero, that $\mu_i=\mu_i'$ for $i<k$ and (2.5) we may again show analogously that the map $f_b^k:\G_-\rightarrow\G$ given by
$$f_b^k(x):=\mu_k(b,x)-\mu'_k(b,x),\quad x\in\G_-$$
defines a Spencer $1$-cocycle.  In a completely analogous fashion we define the map $\rho_1^k:\G_1\rightarrow\G$ such that
\begin{eqnarray}
[\rho_1^k(b),x]'_{\epsilon}&=&[b,x]+\sum_{i<k}\mu'_i(b,x){\epsilon}^i
\nonumber \\
&+&\mu_k(b,x)\epsilon^k+\sum_{i>k}\mu''_i(b,x){\epsilon}^i,\quad\forall b\in\G_1,\forall x\in\G_-.\nonumber
\end{eqnarray}

Now $\rho_j^k$, for $j\ge 2$, are defined analogously. The sequence above, call it $\{\mu''_1,\mu''_2,\cdots\}$, of course is a defining sequence for $\bar{\G}'$.  We have $\mu_i=\mu''_i$ for $i<k$ and $\mu_k|_{\G_-\times\G}=\mu''_k|_{\G_-\times\G}$.  App
lying Proposition 2.5 again we have $\mu_k=\mu_k''$. Also since the subspace $V_-$ remains unchanged, obviously $\mu''_i|_{\G_-\times\G_-}=\mu'_i|_{\G_-\times\G_-}$.  {$\qed$\medskip}

The following is an important remark.

\remark{Remark 2.3.} In the proof of Proposition 2.6 the only place where full prolongation is used is to find elements $v_a$, which then allows us to define $\rho_i^k$.  Now we may assume that $\mu_j|{\A\times\G}=0$ for all $j$, where $\A$ is the maximal
 reductive subalgebra of $(\G_{0})_{\bar{0}}$, as explained earlier.  Using this it is easy to verify that $\rho_i^k$ is an injective $\A$-homomorphism.  Hence the map $a\rightarrow v_a$ is an $\A$-homomorphism.  In particular, if $H^{l,1}(\G_-;\G)$ is a 
direct sum of irreducible $\A$-modules that are not isomorphic to those irreducible $\A$-modules that appear in the decomposition of $\G_{(k)}$, then we may always find such $v_a$'s.  Therefore the assumption of full prolongation of degree $k$ may be repl
aced by the weaker assumption of
$${\rm Hom}_{\A}(H^{l,1}(\G_-;\G),\G_{(1)})=0,\quad\forall l\ge k,$$
and the conclusion of Proposition 2.5 remains valid. We will say that $\G$ is an {\it almost full prolongation} of $\G_{\le 0}$ if
${\rm Hom}_{\A}(H^{l,1}(\G_-;\G),\G_{(1)})=0$ for all $l\ge 1$.

Combining Remark 2.3 with Proposition 2.6 we have proved

\proclaim{Theorem 2.1.} Let $\G=\oplus_{j=-h}^{\infty}\G_j$ be a transitive graded Lie superalgebra.  Let $\bar{\G}_{\epsilon}$ and $\bar{\G}_{\epsilon}'$ be two filtered deformations of $\G$ with defining sequences $\{\mu_1,\mu_2,\cdots\}$ and $\{\mu'_1,
\mu'_2,\cdots\}$, respectively.  Suppose that $\mu_i|_{\G_-\times\G_-}=\mu'_i|_{\G_-\times\G_-}$, for all $i\ge 1$. If furthermore $\G$ is an almost full prolongation of $\G_{\le 0}$, then $\bar{\G}_{\epsilon}\cong\bar{\G}'_{\epsilon}$.

The next two corollaries generalize two results of Kobayashi and Nagano [KN].

\proclaim{Corollary 2.3.} Let $\G=\oplus_{j=-h}^{\infty}\G_j$ be a transitive graded Lie superalgebra.  Suppose that $H^{l,2}(\G_-;\G)_{\bar{0}}$ contains no non-trivial $\G_0$-invariant vectors and that $\G$ is an almost full prolongation of $\G_{\le 0}$
.  Then $\G$ has no non-trivial filtered deformations.

\proof Let $\bar{\G}_{\epsilon}$ correspond to $\{\mu_1,\mu_2,\cdots\}$.  By Propositions 2.2 and 2.3 we may assume that $\mu_i|_{\G_-\times\G_-}=0$.  But then Theorem 2.1 tells us that $\bar{\G}_{\epsilon}$ is isomorphic to the trivial deformation. {$\qe
d$\medskip} 

\proclaim{Corollary 2.4.} Let $L=L_{-h}\supset\cdots\supset L_{-1}\supset L_{0}\supset L_1\supset\cdots$ be a filtered deformation of a transitive graded Lie superalgebra $\G=\oplus_{j=-h}^{\infty}\G_j$.  Suppose that $\G$ is an almost full prolongation o
f $\G_{\le 0}$.  If there exists a $\Z$-graded  subalgebra $V_{-}=\oplus_{j< 0}V_j$ of $L$ isomorphic to $\G_-$ such that $V_{-}+L_0=L$, then $L\cong\bar{\G}$.

\proof The existence of the subalgebra $V_{-}$ means that there exists a defining sequence $\{\mu_1,\mu_2,\cdots\}$ of $L$ such that $\mu_i|_{\G_{-}\times\G_{-}}=0$ for all $i$.  But since it is an almost full prolongation, Theorem 2.1 tells us that $\mu_
{i}=0$ for all $i$.  {$\qed$\medskip}

\proclaim{Corollary 2.5.} Let $\G$ be a transitive graded Lie superalgebra.  Let $\bar{\G}_{\epsilon}$ and $\bar{\G}'_{\epsilon}$ be two filtered deformations given by $\{\mu_1,\mu_2,\cdots\}$ and $\{\mu'_1,\mu'_2,\cdots\}$, respectively. Suppose that $\mu_i|_{\G_-\times\G_-}=\mu'_i|_{\G_-\times\G_-}$ for $i=1,\cdots,k-1$ and $H^{l,2}(\G_-;\G)_{\bar{0}}$ contains no $\G_0$-invariant vectors for $l\ge k$.  If furthermore $\G$ is an almost full prolongation of $\G_{\le 0}$, then $\bar{\G}_{\epsilon}\cong\bar{\G}'_{\epsilon}$.

\proof By Proposition 2.6 we may assume that $\mu_i=\mu_i'$ for $i=1,\cdots,k-1$.  By formulas (2.5) and (2.6) it follows that $(\mu_k-\mu_k')|_{\G_-\times\G_-}$ is a $\G_0$-invariant Spencer cocycle.  Hence by Proposition 2.6 again we may assume that $\mu_k=\mu_k'$. Proceeding this way we show that $\mu_i=\mu_i'$ for all $i\ge k$. {$\qed$\medskip}

\proclaim{Proposition 2.7.} {\rm [G]} Let $\G$ be a transitive $\Z$-graded Lie superalgebra such that $\G_{-2}$ contains a central element of $\G$ of parity $\delta$.  If $H^{1}(\G_0;\G_{-1})_{\delta}=0$, then $\G$ has no filtered deformation $L$ such tha
t $L_{0}$ is a maximal subalgebra.

\proof Let $1$ denote this central element.  Suppose that $\{\mu_1,\mu_2\cdots\}$ is a defining sequence of a deformation $\bar{\G}_{\epsilon}=L$ of $\G$.  We will show that there exists a defining sequence $\{\mu_1',\mu_2'\cdots\}$ of $L$ such that $\mu_1'(1,\G_0)=0$.  From this it follows that $1$ normalizes ${\G}_{(0)}$ and hence $L_{0}$ is not maximal.

For $a,b\in\G_0$ we have by Jacobi identity
$$[1,[a,b]_{\epsilon}]_{\epsilon}=[[1,a]_{\epsilon},b]_{\epsilon}+p(1,a)[a,[1,b]_{\epsilon}]_{\epsilon}.$$
Collecting the coefficient of $\epsilon$ we get
$$\mu_1(1,[a,b])=[\mu_1(1,a),b]+p(1,a)[a,\mu_1(1,b)].$$
This means precisely that the map $c:\G_0\rightarrow\G_{-1}$ given by $c(a)=\mu_1(1,a)$ is a $1$-cocycle of $\G_0$ with coefficients in $\G_{-1}$.  By assumption $c$ is a coboundary, and hence we may add to $1$ an element $x$ in $\G_{-1}$ so that
$$[1+x,a]_{\epsilon}=\mu_2'(1,a)\epsilon^2+o(\epsilon^3).$$
This choice of $V_{-2}$ gives the required defining sequence. {$\qed$\medskip}

\proclaim{Lemma 2.1.} Let $\hat{\G}=\G+\C 1$ be a transitive $\Z$-graded Lie superalgebra, which is a central extension of a $\Z$-graded Lie superalgebra by adding a central element $1$ in degree $-2$ such that $[\G_{-1},\G_{-1}]=\C 1$.  Then $\hat{\G}$ i
s the full prolongation of $\hat{\G}_{\le 0}$, provided that $\G$ is the full prolongation of $\G_{\le 0}$.

\proof  Let $\alpha:\G_{-1}+\C 1\rightarrow\hat{\G}$ be a derivation of degree $\ge 1$.  Thus  $\alpha:\G_{-1}+\C 1\rightarrow\G$.  It suffices to show that $\alpha(1)=0$.  Let $x\in\G_{-1}$.  Then $0=\alpha([x,1])=[\alpha(x),1]+p(\alpha,x)[x,\alpha(1)]=p
(\alpha,x)[x,\alpha(1)]$.  Thus by transitivity $\alpha(1)=\lambda 1$, $\lambda\in\C$.  But then $\lambda=0$, since $\alpha$ is of positive degree. {$\qed$\medskip}

\heading{3}{Examples of $\Z$-graded Lie superalgebras}

In this section we will recall the definitions and list some properties of those $\Z$-graded Lie superalgebras whose filtered deformations we are interested in.  Some of their properties are well-known and can be found in [S].

Let $\Lambda(n)$ be the Grassmann superalgebra in the $n$ odd indeterminates $\xi_1,\xi_2,\cdots,\xi_n$.  Let $x_1,x_2,\cdots,x_m$ be $m$ even indeterminates.  Set $\Lambda(m,n)=\C[x_1,\cdots,x_m]\otimes\Lambda(n)$.  Then $\Lambda(m,n)$ is an associative 
commutative superalgebra.  Let $W(m,n)$ be the Lie superalgebra of derivations of $\Lambda(m,n)$.  Then $W(m,n)$ consists of elements of the form [K1]:
$$\sum_{i=1}^mf_i{{\dd}\over{\dd x_i}}+\sum_{i=1}^n g_i{{\dd}\over{\dd\xi_i}},$$
where $f_i,g_i\in\Lambda(m,n)$ and ${{\partial}\over{\partial x_i}}$ (respectively ${{\partial}\over{\partial\xi_i}}$) is the even (respectively odd) derivation uniquely determined by ${{\partial}\over{\partial x_i}}(x_j)=\delta_{ij}$ and ${{\partial}\over{\partial x_i}}(\xi_j)=0$ (respectively ${{\partial}\over{\partial\xi_i}}(x_j)=0$ and ${{\partial}\over{\partial\xi_i}}(\xi_j)=\delta_{ij}$).  To each vector field $D=\sum_{i=1}^mf_i{{\partial}\over{\partial x_i}}+\sum_{i=1}^n g_i{{\partial}\over{\partial\xi_i}}$ we may associate its divergence by
$${\rm div}D=\sum_{i=1}^m {{\dd f_i}\over{\dd x_i}}+\sum_{i=1}^n(-1)^{p(g_i)}{{\dd g_i}\over{\dd\xi_i}}.$$

Let $\Omega(m,n)$ be the superalgebra of differential forms over $\Lambda(m,n)$ [K1].  Consider the following differential form:
$$\omega=\sum_{i=1}^n dx_i d\xi_i\in\Omega(n,n).$$  
Define the odd Hamiltonian superalgebra [L]
$$HO(n,n):=\{D\in W(n,n)|D\omega=0\}.$$
The Lie superalgebra $HO(n,n)$ is simple if and only if $n\ge 2$.

The Lie superalgebra $HO(n,n)$ contains the subalgebra of divergence free vector fields
$$SHO'(n,n):=\{D\in HO(n,n)|{\rm div}D=0\}.$$
The derived algebra of $SHO'(n,n)$ is an ideal of codimension $1$, denoted by $SHO(n,n)$, provided that $n\ge 2$.  $SHO(n,n)$ is simple if and only if $n\ge 3$.

In $\Lambda(n,n)$ we can define the Buttin bracket by
$$[f,g]:=\sum_{i=1}^n({{\dd f}\over{\dd x_i}}{{\dd
g}\over{\dd\xi_i}}+(-1)^{p(f)}{{\dd f}\over{\dd\xi_i}}{{\dd
g}\over{\dd x_i}}),$$ which makes $\Lambda(n,n)$, with reversed parity, into a Lie superalgebra.  It contains a one-dimensional center consisting of
constant functions.  The map $\Lambda(n,n)\rightarrow HO(n,n)$ given
by $$f\rightarrow \sum_{i=1}^n({{\dd f}\over{\dd
x_i}}{{\dd}\over{d\xi_i}}+(-1)^{p(f)}{{\dd
f}\over{\dd\xi_i}}{{\dd}\over{\dd x_i}}),$$ is a surjective
homomorphism of Lie superalgebras with kernel consisting of constant
functions.  Hence we may (and will) identify $HO(n,n)$ with
$\Lambda(n,n)/\C 1$ with reversed parity.  In this identification we have:
$$SHO'(n,n)=\{f\in\Lambda(n,n)/\C 1|\ \Delta(f)=0\},$$ where
$\Delta=\sum_{i=1}^n{{\partial^2}\over{\partial x_i\partial\xi_i}}$ is
the odd Laplace operator, and $SHO(n,n)$ is identified with the
subspace consisting of elements not containing the monomial
$\xi_1\xi_2\cdots\xi_n$.

By putting ${\deg}x_i=1$ and ${\rm deg}\xi_i=1$ for $i=1,2,\cdots,n$ the Lie superalgebras $HO(n,n)$, $SHO'(n,n)$ and $SHO(n,n)$ become $\Z$-graded Lie superalgebras of depth $1$.  Since $[x_i,\xi_j]=\delta_{ij}1$ we obtain, by adding $\C 1$ to degree $-2
$, non-trivial central extensions of $HO(n,n)$, $SHO'(n,n)$ and $SHO(n,n)$, denoted by $\hat{HO}(n,n)$, $\hat{SHO}'(n,n)$ and $\hat{SHO}(n,n)$, respectively.

The $0$-th graded components of $HO(n,n)$ and $\hat{HO}(n,n)$ have a basis consisting of vectors of the form $\{x_ix_j\}$, $\{x_i\xi_j\}$ and $\{\xi_i\xi_j\}_{i\not=j}$ for $i,j=1,2,\cdots,n$.  This is the $\Z$-graded finite-dimensional Lie superalgebra $
\tilde{P}(n)=\tilde{P}(n)_{-1}+\tilde{P}(n)_0+\tilde{P}(n)_1$ (cf.~[K1]), where $\tilde{P}(n)_{0}\cong gl_n$, $\tilde{P}(n)_{-1}\cong\Lambda^2(\C^{n*})$ and $\tilde{P}(n)_{1}\cong S^2(\C^n)$, where $\C^n$ stands for the standard representation of $gl_n$. 
 Their $-1$-st graded components have a basis consisting of $\{x_i\}$ and $\{\xi_i\}$, $i=1,2,\cdots,n$.  Evidently the span of $\{x_i\}$ as a $gl_n$-module is isomorphic to $\C^n$, while the span of $\{\xi_i\}$ is isomorphic to $\C^{n*}$.

The $0$-th graded components of $SHO(n,n)$, $SHO'(n,n)$, $\hat{SHO}(n,n)$ and $\hat{SHO}'(n,n)$ have a basis consisting of vectors of the form $\{x_ix_j\}$, $\{x_i\xi_j\}_{i\not=j}$, $\{x_i\xi_i-x_{i+1}\xi_{i+1}\}_{i<n}$ and $\{\xi_i\xi_j\}_{i\not=j}$ for
 $i,j=1,2,\cdots,n$.  This is the $\Z$-graded subalgebra $P(n)$ of $\tilde{P}(n)$ (cf.~[K1]), where $P(n)_{0}\cong sl_n$ and $P(n)_{-1}\cong\Lambda^2(\C^{n*})$ and $P(n)_{1}\cong S^2(\C^n)$, where $\C^n$ stands for the standard representation of $sl_n$.  
Similarly, their $-1$-st graded components have a basis consisting of $\{x_i\}$ and $\{\xi_i\}$, $i=1,2,\cdots,n$ with the span of $\{x_i\}$ isomorphic to $\C^n$ and the span of $\{\xi_i\}$ isomorphic to $\C^{n*}$.

It can be shown that $HO(n,n)$ and $SHO'(n,n)$ are full prolongations of $\sum_{i=1}^n(\C x_i+\C\xi_i)\oplus\tilde{P}(n)$ and $\sum_{i=1}^n(\C x_i+\C\xi_i)\oplus P(n)$, respectively [CK].  Thus it follows from Lemma 2.1 that $\hat{HO}(n,n)$ and $\hat{SHO}
'(n,n)$ are full prolongations of $\C 1\oplus\sum_{i=1}^n(\C x_i+\C\xi_i)\oplus\tilde{P}(n)$ and $\C 1\oplus\sum_{i=1}^n(\C x_i+\C\xi_i)\oplus P(n)$, respectively, as well.  Consequently, in the case of $SHO(n,n)$ and $\hat{SHO}(n,n)$ the only obstruction
 to being a full prolongation lies in degree $n-2$.  More precisely we have ($l\in\Z_+$)
\begin{eqnarray}
&&H^{l,1}(SHO(n,n)_{-1};SHO(n,n))=0,\quad l\not=n-2,\nonumber \\
&&H^{l,1}(\hat{SHO}(n,n)_{-};\hat{SHO}(n,n)) =0,\quad l\not=n-2,\nonumber \\
&&H^{n-2,1}(SHO(n,n)_{-1};SHO(n,n))=\C\xi_1\xi_2\cdots\xi_n,\nonumber \\
&&H^{n-2,1}(\hat{SHO}(n,n)_{-};\hat{SHO}(n,n)) =\C\xi_1\xi_2\cdots\xi_n. \;\;\;\;\;\;\;\;\;\;\;\;\;\;\;\;\;\;\;\;\;\;\;\;\;\;\;\; (3.1)\nonumber
\end{eqnarray}
Let us denote the vector $\sum_{i=1}^n x_i\xi_i$ in $HO(n,n)_0$ by $\Phi$.  Other Lie superalgebras that arise in the classification in [K3], and hence whose filtered deformations we also need to consider are the following four series: $SHO(n,n)+\C\Phi$, 
$SHO'(n,n)+\C\Phi$, $\hat{SHO}(n,n)+\C\Phi$ and $\hat{SHO}'(n,n)+\C\Phi$.

Let $x_1,x_2,\cdots,x_n$ be $n$ even indeterminates and $\xi_1,\xi_2,\cdots,\xi_n,\xi_{n+1}=\tau$ be $n+1$ odd indeterminates.  Define the odd contact form to be
$$\Omega=d\tau+\sum_{i=1}^n(\xi_idx_i+x_id\xi_i)\in\Omega(n,n+1).$$
The odd contact superalgebra $KO(n,n+1)$ is the following subalgebra of $W(n,n+1)$ [ALS]:
$$KO(n,n+1)=\{D\in W(n,n+1)| D\Omega=f_{D}\Omega\},$$
for some $f_D\in\C[\tau,x_1,\cdots,x_n,\xi_1,\cdots,\xi_n]$.    The Lie superalgebra $KO(n,n+1)$ can be realized as follows.  We may define the odd contact bracket on the space $\Lambda(n,n+1)$ by
$$[f,g]=(2-E)f{{\partial g}\over{\partial\tau}}+(-1)^{p(f)}{{\partial f}\over{\partial\tau}}(2-E)g-\sum_{i=1}^n({{\dd f}\over{\dd x_i}}{{\dd g}\over{\dd\xi_i}}+(-1)^{p(f)}{{\dd f}\over{\dd\xi_i}}{{\dd g}\over{\dd x_i}}),$$
where $E=\sum_{i=1}^n({{\partial}\over {\partial x_i}}+{{\partial}\over {\partial \xi_i}})$ is the Euler operator.  Reversing parity, $\Lambda(n,n+1)$ with this bracket becomes a Lie superalgebra and the map $\Lambda(n,n+1)\rightarrow KO(n,n+1)$, given by

$$f\rightarrow (2-E)f{{\partial}\over{\partial\tau}}-(-1)^{p(f)}{{\partial f}\over{\partial\tau}}E-\sum_{i=1}^n({{\dd f}\over{\dd x_i}}{{\dd}\over{d\xi_i}}+(-1)^{p(f)}{{\dd f}\over{\dd\xi_i}}{{\dd}\over{\dd x_i}}),$$
is a isomorphism of Lie superalgebras.  Hence we may (and will) identify the Lie superalgebra $KO(n,n+1)$ with $\Lambda(n,n+1)$ with reversed parity.

For $\beta\in\C$ we let ${\rm div}_{\beta}=2(-1)^{p(f)}(\Delta+(E-n\beta){{\partial}\over{\partial\tau}})$, where $\Delta$ is the odd Laplace operator.  We set [Ko]
$$SKO'(n,n+1;{\beta})=\{f\in\Lambda(n,n+1)|\ {\rm div}_{\beta}f=0\}.$$
This is a subalgebra of $KO(n,n+1)$ and is simple if and only if $n\ge 2$ and 
$\beta\not=1,{{n-2}\over n}$.  Let $SKO(n,n+1;\beta)$ denote the derived algebra of $SKO'(n,n+1;\beta)$.  Then $SKO(n,n+1;\beta)$ is simple for $n\ge 2$
and it coincides with $SKO'(n,n+1;\beta)$ unless $\beta=1$ or $\beta={{n-2}\over n}$.  The Lie superalgebra $SKO(n,n+1;1)$ (respectively $SKO(n,n+1;{{n-2}\over n})$) consists of elements not containing the monomial $\tau\xi_1\xi_2\cdots\xi_n$ (respectivel
y $\xi_1\xi_2\cdots\xi_n$).
Note that $SKO(n,n+1;{1\over n})$ is the subalgebra of $KO(n,n+1)$ consisting of divergence free vector fields.

By setting ${\rm deg}\tau=2$ and ${\rm deg}x_i={\rm deg}\xi_i=1$ for all $i$, $KO(n,n+1)$ and hence $SKO(n,n+1;\beta)$ and $SKO'(n,n+1;\beta)$ (since ${\rm div}_{\beta}$ is homogeneous with respect to this gradation) become $\Z$-graded Lie superalgebras. 
 They are all of depth $2$.  In the case of $KO(n,n+1)$ the $0$-th graded component has a non-trivial center, namely $\C\tau$, and hence by Corollary 2.2, $KO(n,n+1)$ has no non-trivial filtered deformations.

Now consider $SKO(n,n+1;\beta)$ and $SKO'(n,n+1;\beta)$.  The $0$-th graded component is spanned by the vectors $\{x_ix_j\}$, $\{x_i\xi_j\}$, $\{\xi_i\xi_j\}_{i\not=j}$ and $\tau+\beta\Phi$, where $i,j=1,2,\cdots,n$ and $\Phi=\sum_{i=1}^nx_i\xi_i$.  This 
is the Lie superalgebra $\tilde{P}(n)=P(n)+\C(\tau+\beta\Phi)$.  The $-2$-nd graded component is spanned by $\C 1$, on which $\tau+\beta\Phi$ acts as the scalar $-2$.  The $-1$-st graded component is spanned by the vectors $\{x_i\}$ and $\{\xi_i\}$ for $i
=1,\cdots,n$.  With respect to $P(n)$ this is the standard representation, and $\tau+\beta\Phi$ acts on $\sum_{i=1}^n\C x_i$ (respectively $\sum_{i=1}^n\C \xi_i$) as the scalar $-1+\beta$ (respectively $-1-\beta$).  It can be shown that $SKO'(n,n+1;\beta)
$ are full prolongations for all $\beta$ [CK].  Hence we have:
\newpage
\begin{eqnarray}
&&H^{l,1}(SKO(n,n+1;1)_{-1};SKO(n,n+1;1))=0,\quad l\not=n,\quad \quad \;\;\;\;\;\;(3.2)\nonumber \\
&&H^{l,1}(SKO(n,n+1;{{n-2}\over n})_{-1};SKO(n,n+1;{{n-2}\over n}))=0,\quad \nonumber l\not=n-2,\nonumber \\
&&H^{n,1}(SKO(n,n+1;1)_{-1};SKO(n,n+1;1))=\C\tau\xi_1\xi_2\cdots\xi_n,\nonumber \\
&&H^{n-2,1}(SKO(n,n+1;{{n-2}\over n})_{-1};SKO(n,n+1;{{n-2}\over n}))=\C\xi_1\xi_2\cdots\xi_n.\nonumber \end{eqnarray}

Let $p_1,p_2,\cdots,p_n,q_1,q_2,\cdots,q_n$ be $2n$ even and $\xi_1,\xi_2,\cdots,\xi_m$ be $m$ odd variables.  Consider the differential form
$$\sigma=\sum_{i=1}^ndp_i dq_i+\sum_{i=1}^md\xi_i d\xi_i\in\Omega(2n,m).$$
Define the Hamiltonian superalgebra to be [K1]
$$H(2n,m)=\{D\in W(2n,m)|D\sigma=0\}.$$
It is a simple Lie superalgebra for $n\ge 1$ and $m\ge 0$.

Let $\Lambda(2n,m)=\C[p_1,\cdots,p_n,q_1,\cdots,q_n]\otimes\Lambda(m)$.  For $f,g\in\Lambda(2n,m)$ we define the Poisson bracket
$$[f,g]=\sum_{i=1}^n({{\partial f}\over{\partial p_i}}{{\partial g}\over{\partial q_i}}-{{\partial f}\over{\partial q_i}}{{\partial g}\over{\partial p_i}})-(-1)^{p(f)}\sum_{i=1}^m{{\partial f}\over{\partial \xi_i}}{{\partial g}\over{\partial \xi_i}}.$$
As before $\Lambda(2n,m)$ with this Poisson bracket is a Lie superalgebra.  The map
$$f\rightarrow\sum_{i=1}^n({{\partial f}\over{\partial p_i}}{{\partial}\over{\partial q_i}}-{{\partial f}\over{\partial q_i}}{{\partial}\over{\partial p_i}})-(-1)^{p(f)}\sum_{i=1}^m{{\partial f}\over{\partial \xi_i}}{{\partial}\over{\partial \xi_i}}$$
defines a surjective homomorphism of Lie superalgebras from $\Lambda(2n,m)$ onto $H(2n,m)$.  The kernel of this map consists of constant functions so that we may (and will) identify $H(2n,m)$ with $\Lambda(2n,m)/\C 1$.  By putting ${\rm deg}p_i={\rm deg}q
_i=1$ and ${\rm deg}\xi_j=1$ for $i=1,\cdots,n$ and $j=1,\cdots,m$ $H(2n,m)$ becomes a $\Z$-graded Lie superalgebra of depth $1$.

 The $0$-th graded component of $H(2n,m)$ is the Lie superalgebra $spo(2n,m)$.  Now $spo(2n,m)_{\bar{0}}$ has a basis consisting of vectors of the form $\{p_ip_j,p_iq_j,q_iq_j\}$ for ${i,j=1,\cdots,n}$ and $\{\xi_i\xi_j\}_{i\not=j}$ for ${i,j=1,\cdots,m}$
 and hence is isomorphic to the Lie algebra $sp_{2n}\oplus so_m$.  $spo(2n,m)_{\bar{1}}$ has a basis consisting of vectors of the form $\{p_i\xi_j,q_i\xi_j\}$ for ${i=1,\cdots,n}$ and ${j=1,\cdots,m}$.  Its span is isomorphic to the $sp_{2n}\oplus so_m$-m
odule $\C^{2n}\otimes \C^m$, where $\C^{2n}$ and $\C^m$ are the respective standard representations of $sp_{2n}$ and $so_m$.  $H(2n,m)_{-1}$ has a basis consisting of vectors of the form $\{p_i,q_i\}$ and $\{\xi_j\}$, $i=1,\cdots,n$ and $j=1,\cdots,m$.  E
vidently the span of $\{p_i,q_i\}$ is isomorphic to $\C^{2n}$, while the span of $\{\xi_j\}$ is isomorphic to $\C^m$.  It is the standard representation of $spo(2n,m)$, denoted by $\C^{2n|m}$.

Finally $H(2n,m)$ is the full prolongation of the pair $\C^{2n|m}\oplus spo(2n,m)$ [CK].
\heading{4}{Calculations of Spencer $2$-cocycles and triviality of filtered deformations of $SHO(n,n)$, $HO(n,n)$, $H(2m,n)$ and $SKO(n,n+1;\beta)$, $\beta\not={{n+2}\over n}$}

In this section we will apply the results obtained in Section 2 and start our investigations of filtered deformations of those graded Lie superalgebras discussed in Section 3.  Due to Proposition 2.2 our first step should be to find $\G_0$-invariant Spenc
er $2$-cocycles.  However, because of lack of complete reducibility of $\G_0$-modules in general, we will restrict ourselves to the even part $\A=(\G_0)_{\bar{0}}$ of $\G_0$, for which we do have complete reducibility (in all our examples).  So our task w
ill be first to look for all $\A$-invariant vectors on the level of $2$-cochains and then determine which of those vectors indeed satisfy the $2$-cocycle condition.  To find $\A$-invariant $2$-cochains we will first need the $\A$-module structure of $\G_j
$ and then use this to find all the trivial $\A$-modules that appear in $\Lambda^2(\G_-^*)\otimes\G_j$.  So our first task will be to determine the $\A$-module structure of $\G_j$ for every $j$.

Consider the Lie superalgebra $\G=SHO'(n,n)$, for $n\ge 3$.  As usual we will write $\G_j$ for its $j$-th graded component.  Here $\A=sl_n$ and we denote by $R(\sum_{i}k_i\pi_i)$ the irreducible $sl_n$-module with highest weight $\sum_{i}k_i\pi_i$, where 
$\pi_i$ is the $i$-th fundamental weight of $sl_n$.  Below we will list the structure of $\G_j$ as $sl_n$-modules and also include explicitly a highest weight vector.

\begin{eqnarray}
\G_{-1}&&:\{R(\pi_1), x_1\},\{R(\pi_{n-1}),\xi_n\}\nonumber \\
\G_{0}&&:\{R(2\pi_1),x_1^2\},\nonumber \\
&&\{R(\pi_1+\pi_{n-1}),x_1\xi_n\},\nonumber \\
&&\{R(\pi_{n-2}),\xi_n\xi_{n-1}\}\nonumber \\
\G_{1}&&:\{R(3\pi_1),x_1^3\},\{R(2\pi_1+\pi_{n-1}),x_1^2\xi_n\},\nonumber \\
&&\{R(\pi_1+\pi_{n-2}),x_1\xi_{n}\xi_{n-1}\},\nonumber \\
&&\{R(\pi_{n-3}),\xi_n\xi_{n-1}\xi_{n-2}\}\nonumber \\
&&\qquad\vdots\nonumber \\
\G_{n-2}&&:\{R(n\pi_1),x_1^n\},\{R((n-1)\pi_1+\pi_{n-1}),x_1^{n-1}\xi_n\},\nonumber \\
&&\{R((n-2)\pi_1+\pi_{n-2}),x_1^{n-2}\xi_{n}\xi_{n-1}\},\cdots,\nonumber \\
&&\{R(2\pi_1),x_1\xi_n\cdots\xi_2\},\{R(0),\xi_n\xi_{n-1}\cdots\xi_{1}\}\nonumber \\ \G_{n-1}&&:\{R((n+1)\pi_1),x_1^{n+1}\},\{R(n\pi_1+\pi_{n-1}),x_1^{n}\xi_n\},\nonumber \\ &&\{R((n-1)\pi_1+\pi_{n-2}),x_1^{n-1}\xi_{n}\xi_{n-1}\},\cdots,\nonumber \\
&&\{R(3\pi_1),x_1^2\xi_n\cdots\xi_2\}\nonumber \\
&&\qquad\vdots\nonumber 
\end{eqnarray}
Note that the structure as an $sl_n$-module of $SHO(n,n)$ is exactly the same except that the component $\{R(0),\xi_n\xi_{n-1}\cdots\xi_{1}\}$ is removed in $\G_{n-2}$.

For $i=1,\cdots,n$ let $f_i,\theta_i\in\G_{-1}^*$ be such that $f_i(x_j)=\delta_{ij}$, $f_i(\xi_j)=0$, $\theta_i(x_j)=0$ and $\theta_i(\xi_j)=\delta_{ij}$.  Note that since the map $\Lambda(n,n)\rightarrow HO(n,n)$ is odd (see Section 3), we have to rever
se parity in $\Lambda(n,n)$ in order to get the correct parity. Hence $p(f_i)=\bar{1}$ and $p(\theta_i)=\bar{0}$.  Evidently the span of $\{f_i\}$ is the $sl_n$-module $R(\pi_{n-1})$ with highest weight vector $f_n$, while the span of $\{\theta_i\}$ is $R
(\pi_1)$ with highest weight vector $\theta_1$.  So using our notation
$\Lambda^2(\G^*_{-1})$ consists of the following irreducible components with highest weight vectors:
$$\{R(2\pi_{n-1}),f_n^2\},\{R(\pi_1+\pi_{n-1}),\theta_1 f_n\},\{R(\pi_{2}),\theta_1\theta_2\},\{R(0),\sum_{i=1}^nf_i\theta_i\}.$$
So to find trivial $sl_n$-modules in $\Lambda^2(\G_{-1}^*)\otimes\G_{j}$ we need to find modules in the table above that are contragredient to these modules in $\Lambda^2(\G_{-1}^*)$.

Let $\G=SHO(n,n)$ or $\G=SHO'(n,n)$ $n\ge 3$.
In the case when $n=3$ $\G_{-1}\cong R(\pi_1)\oplus R(\pi_2)$.  Hence there exits a trivial $sl_3$-module in $\Lambda^2(\G_{-1}^*)\otimes\G_{-1}$, given by the vector $\sum_{i=1}^3\theta^*_i\otimes x_i$, where $\theta_i^*$ stands for the Hodge dual of $\theta_i$.  However this vector is odd, hence it cannot give a deformation.  For $n>3$ there are no trivial $sl_n$-modules in $\Lambda^2(\G_{-1}^*)\otimes\G_{-1}$.

Now for $n\ge 3$ we have the following linearly independent $sl_n$-invariant vectors in $\Lambda^2(\G_{-1}^*)\otimes\G_0$:
\begin{eqnarray}
c_1&=&\sum_{i,j}f_i f_j\otimes x_i x_j,\nonumber \\
c_2&=&{1\over 2}\sum_{i=1}^{n-1}(f_i\theta_i-f_{i+1}\theta_{i+1})\otimes(x_i\xi_i-x_{i+1}\xi_{i+1})\nonumber \\&+&\sum_{i\not=j}f_i\theta_j\otimes x_i\xi_j,\nonumber \\
c_3&=&\sum_{i<j}\theta_i\theta_j\otimes\xi_i\xi_j.\nonumber 
\end{eqnarray}

But $c_1,c_2,c_3$ are all odd vectors, hence they cannot give deformations.

In $\Lambda^2(\G_{-1}^*)\otimes\G_{n-2}$ we have the following $sl_n$-invariant vector in $SHO(n,n)$:
$$c_1=\sum_{i,j}f_i f_j\otimes x_i\xi^*_j,$$
where $\xi_j^*$ again denotes the Hodge dual of $\xi_j$.  However,
$$dc_1(\xi_1,x_1,x_1)=-[\xi_1,x_1\xi_1^*]\not=0.$$
Hence $c_1$ is not a cocycle.  In the case of $SHO'(n,n)$ there is an additional $sl_n$-invariant vector in $\Lambda^2(\G_{-1}^*)\otimes\G_{n-2}$ that is not proportional to $c_1$, namely
$$c_2=\sum_{i=1}^n f_i\theta_i\otimes\xi_1\xi_2\cdots\xi_n.$$
It is subject to direct verification that $c_1+c_2$ is an even Spencer cocycle for $n$ even.  For $n$ odd, obviously $c_1$ and $c_2$ are both odd.  We summarize the above computation.

\proclaim{Proposition 4.1.}
\begin{eqnarray}
&&H^{l,2}(SHO(n,n)_{-1};SHO(n,n))_{\bar{0}}^{sl_n}=0,\quad n\ge 3,\nonumber \\
&&H^{l,2}(SHO'(n,n)_{-1};SHO'(n,n))_{\bar{0}}^{sl_n}=0,\quad n\ge 3,\ l\not=n,\nonumber \\
&&H^{n,2}(SHO'(n,n)_{-1};SHO'(n,n))_{\bar{0}}^{sl_n}=\C,\quad n\ge 3,\ n\ {\rm even},\nonumber \\
&&H^{n,2}(SHO'(n,n)_{-1};SHO'(n,n))_{\bar{0}}^{sl_n}=0,\quad n\ge 3,\ n\ {\rm odd.}\nonumber 
\end{eqnarray}

From the structure of $SHO(n,n)$ as an $sl_n$-module given in the above table, we see that the trivial $sl_n$-module doesn't appear in the decomposition of $SHO(n,n)$, for $n\ge 3$.  Hence by (3.1) $SHO(n,n)$ is an almost full prolongation and so by Corol
lary 2.3 we obtain

\proclaim{Theorem 4.1.} $SHO(n,n)$ has no non-trivial filtered deformations for $n\ge 3$.

Next consider the case of $SHO'(2,2)$.  In this case the maximal reductive subalgebra of $\G_0$ is $sl_2$.  Denoting by $R(k)$ the irreducible $sl_2$-module of highest weight $k\in\Z_+$, $SHO'(2,2)$ decomposes as an $sl_2$-module as follows:
\begin{eqnarray}
\G_{-1}&:&\{R(1), x_1\},\{R(1),\xi_2\}\nonumber \\
\G_{0}&:&\{R(2),x_1^2\},\{R(2),x_1\xi_2\},\{R(0),\xi_1\xi_2\}\nonumber \\
\G_{1}&:&\{R(3),x_1^3\},\{R(3),x_1^2\xi_2\}\nonumber \\
&&\qquad\vdots\nonumber \\
\G_{k}&:&\{R((k+2)),x_1^{k+2}\},\{R((k+2)),x_1^{k+1}\xi_2\}\nonumber \\
&&\qquad\vdots\nonumber 
\end{eqnarray}

$\Lambda^2(\G_{-1}^*)\cong\{ R(2),f_2^2\} \oplus\{ R(2),f_2\theta_1\} \oplus\{ R(0),f_1\theta_1+f_2\theta_2\}\oplus\{ R(0),\theta_1\theta_2\}$.  Here as usual $f_i(x_j)=\delta_{ij}$ and $\theta_i(\xi_j)=\delta_{ij}$, $i,j=1,2$.

Thus $R(0)$ can appear in $\Lambda^2(\G_{-1}^*)\otimes\G$ only in degree $2$, i.e.~in $\Lambda^2(\G_{-1}^*)\otimes\G_0$.  The following linearly independent vectors span the even $sl_2$-invariant subspace in $\Lambda^2(\G_{-1}^*)\otimes\G_0$:
\begin{eqnarray}
c_1&=&f_1f_2\otimes(x_1\xi_1-x_2\xi_2)-f_1^2\otimes x_1\xi_2+f_2^2\otimes x_2\xi_1,\nonumber \\
c_2&=&(f_1\theta_1-f_2\theta_2)\otimes x_1x_2-f_1\theta_2\otimes x_1^2+f_2\theta_1\otimes x_2^2,\nonumber \\
c_3&=&(f_1\theta_1+f_2\theta_2)\otimes\xi_1\xi_2.
\nonumber 
\end{eqnarray}

Compute
\begin{eqnarray}
&&dc_1(\xi_1,\xi_2,x_1)=0,\nonumber \\
&&dc_2(\xi_1,\xi_2,x_1)=-[\xi_1,-x_1^2]+[\xi_2,x_1x_2]=-3x_1\not=0,\nonumber \\
&&dc_3(\xi_1,\xi_2,x_1)=[\xi_2,\xi_1\xi_2]=0.\nonumber 
\end{eqnarray}
Hence any cocycle must be a linear combination of $c_1$ and $c_3$.  Now $dc_1(\xi_1,x_1,$\newline$x_1)=[\xi_1,x_1\xi_2]=-\xi_2\not=0$.  One checks easily that $c_1-2c_3$ is a Spencer cocycle. This calculation also shows that $SHO(2,2)$ has no non-trivial 
$sl_2$-invariant Spencer $2$-cocyles. Thus we arrive at

\proclaim{Proposition 4.2.}  Proposition 4.1 holds for $SHO'(2,2)$ and $SHO(2,2)$ as well.

Consider now $\G=\hat{SHO}(n,n)$.  We will now compute $H^{k,2}(\G_-;\G)_{\bar{0}}^{sl_n}$.  Let $c_k\in H^{k,2}(\G_-;\G)^{sl_n}_{\bar{0}}$.  Recalling that $\G_-=\C 1\oplus\G_{-1}$, $c_k$ can be written as
$$c_k=c_k^{-2,-2}+c_k^{-2,-1}+c_k^{-1,-1},$$
where $c_k^{-2,-2}:\C 1\times\C 1\rightarrow\G_{k-4}$, $c_k^{-2,-1}:\C 1\times\G_{-1}\rightarrow\G_{k-3}$ and $c_k^{-1,-1}:\G_{-1}\times\G_{-1}\rightarrow\G_{k-2}$ are bilinear maps.
It follows from the fact that $1$ is central and transitivity that $c^{-2,-2}=0$. Hence
$$c_k=c^{-2,-1}_k+c^{-1,-1}_k.$$
Now both $c^{-2,-1}_k$ and $c^{-1,-1}_k$ are $sl_n$-invariant.  We have investigated $c^{-1,-1}_k$ in our computation Spencer $2$-cocycles of $SHO(n,n)$.  Now we need to investigate $c^{-2,-1}_k$.  From our table of $sl_n$-module structure of $SHO(n,n)$ a
bove, $c^{-2,-1}_k$ can be non-zero only for $k=2$ and $k=n$.  Let $1^*$ denote the dual to $1$.  Then we have the following choices for a non-zero $c^{-2,-1}_k$:

For $k=2$: $c^{-2,-1}_{2}$ is a linear combination of $\sum_{i=1}^n1^*f_i\otimes x_i$ and $\sum_{i=1}^n1^*\theta_i\otimes\xi_i$.  However, both vectors are odd.  So this cannot happen. 

For $k=n$: $c^{-2,-1}_{n}$ is a scalar multiple of
$$c_3=\sum_{i=1}^n1^*f_i\otimes\xi_i^*,$$
where as usual $\xi_i^*$ is the Hodge dual of $\xi_i$.  In this case we see from our previous calculations that $c^{-1,-1}_n$ is a scalar multiple of
$$c_1=\sum_{i,j}f_if_j\otimes x_i\xi_j^*.$$
We know that $c_1$ is not a cocycle.  However, it is easy to verify that $c_1+c_3$ is a non-trivial Spencer cocycle.  Thus we have

\proclaim{Proposition 4.3.} Let $n\ge 2$.  Then 
\begin{eqnarray}
&&H^{l,2}(\hat{SHO}(n,n)_{-};\hat{SHO}(n,n))_{\bar{0}}^{sl_n}=0,\quad l\not=n.\nonumber\\
&&H^{n,2}(\hat{SHO}(n,n)_{-};\hat{SHO}(n,n))_{\bar{0}}^{sl_n}=\C,\quad n\ {\rm even.}\nonumber \\
&&H^{n,2}(\hat{SHO}(n,n)_{-};\hat{SHO}(n,n))_{\bar{0}}^{sl_n}=0,\quad n\ {\rm odd.}\nonumber
\end{eqnarray}

Now consider $\G=\hat{SHO}'(n,n)$.  In this case the calculations of Spencer cocycles is analogous to the case of $\hat{SHO}(n,n)$.  We write an element $c_k\in H^{k,2}(\G_-;\G)$ as $c_k=c^{-2,-1}_k+c^{-1,-1}_k$.  As before $c^{-2,-1}_k$ is non-zero only 
for $k=2,n$.  When $k=2$ this cannot happen as in the case of $\hat{SHO}(n,n)$.  For $k=n$ we conclude as before that $c^{-2,-1}_{n}$ is a scalar multiple of $c_3=\sum_{i=1}^n1^*f_i\otimes\xi_i^*$.  However, there are two other linearly independent $sl_n$
-invariant cochains of degree $n$, as we have seen in the computation of $SHO'(n,n)$, namely $c_1=\sum_{i,j}f_if_j\otimes x_i\xi_j^*$ and $c_2=\sum_{i=1}^nf_i\theta_i\otimes\xi_1\cdots\xi_n$.  There are two linearly independent cocycles in the span of the
se three vectors, namely $c_1+c_2$ and $c_1+c_3$.  However, $c_2-c_3=db$, where $b$ is the Spencer $1$-cochain defined by $b(1)=\xi_1\cdots\xi_n$ and $b(\G_{-1})=0$.

\proclaim{Proposition 4.4.} Let $n\ge 2$. Then
\begin{eqnarray}
&&H^{l,2}(\hat{SHO}'(n,n)_{-};\hat{SHO}'(n,n))_{\bar{0}}^{sl_n}=0,\quad l\not=n.\nonumber \\
&&H^{n,2}(\hat{SHO}'(n,n)_{-};\hat{SHO}'(n,n))_{\bar{0}}^{sl_n}=\C,\quad n\ {\rm even.}\nonumber \\
&&H^{n,2}(\hat{SHO}'(n,n)_{-};\hat{SHO}'(n,n))_{\bar{0}}^{sl_n}=0,\quad n\ {\rm odd.}\nonumber 
\end{eqnarray}

We will now consider $HO(n,n)$.  Here the even part of $\G_0$ is $gl_n$.  First we will assume that $n\not=2,4$.  As an $sl_n$-module $HO(n,n)$ decomposes as follows ($\Phi=\sum_{i=1}^nx_i\xi_i$):
\begin{eqnarray}
\G_{-1}&:&\{R(\pi_1), x_1\},\{R(\pi_{n-1}),\xi_n\}\nonumber \\
\G_{0}&:&\{R(2\pi_1),x_1^2\},\{R(\pi_1+\pi_{n-1}),x_1\xi_n\},
\{R(0),\Phi\},\nonumber \\ &&\{R(\pi_{n-2}),\xi_n\xi_{n-1}\}\nonumber \\
\G_{1}&:&\{R(3\pi_1),x_1^3\},\{R(2\pi_1+\pi_{n-1}),x_1^2\xi_n\},
\{R(\pi_1),x_1\Phi\},\nonumber \\ &&\{R(\pi_1+\pi_{n-2}),x_1\xi_{n}\xi_{n-1}\},\{R(\pi_{n-1}),\xi_n\Phi\},\nonumber \\ &&\{R(\pi_{n-3}),\xi_n\xi_{n-1}\xi_{n-2}\}\nonumber \\
\G_{2}&:&\{R(4\pi_1),x_1^4\},\{R(3\pi_1+\pi_{n-1}),x_1^3\xi_n\},\{R(2\pi_1),x_1^2\Phi\},\nonumber \\ &&\{R(2\pi_1+\pi_{n-2}),x_1^2\xi_{n}\xi_{n-1}\},\{R(\pi_1+\pi_{n-1}),x_1\xi_n\Phi\},\nonumber \\ && \{R(\pi_1+\pi_{n-3}),x_1\xi_n\xi_{n-1}\xi_{n-2}\},\{R(
\pi_{n-2}),\xi_{n-1}\xi_{n-2}\Phi\},\nonumber \\ && \{R(\pi_{n-4}),\xi_n\xi_{n-1}\xi_{n-2}\xi_{n-3}\}\nonumber \\
&&\qquad\vdots\nonumber \\
\G_{n-2}&:&\{R(n\pi_1),x_1^n\},\{R((n-1)\pi_1+\pi_{n-1}),x_1^{n-1}\xi_n\},
\nonumber \\
&&\{R((n-2)\pi_1),x_1^{n-2}\Phi\},\cdots,\nonumber \\
&&\{R(2\pi_1),x_1\xi_n\cdots\xi_2\},\{R(\pi_2),\xi_n\cdots\xi_3\Phi\},\nonumber \\
&&\{R(0),\xi_n\xi_{n-1}\cdots\xi_{1}\}\nonumber \\
\G_{n-1}&:&\{R((n+1)\pi_1),x_1^{n+1}\},\{R(n\pi_1+\pi_{n-1}),x_1^n\xi_n\},\nonumber \\
&&\{R((n-1)\pi_1),x_1^{n-1}\Phi\},\cdots,\{R(3\pi_1),x_1^2\xi_n\cdots\xi_2\},\nonumber \\ &&\{R(\pi_1+\pi_2),x_1\xi_n\cdots\xi_3\Phi\},\{R(\pi_1),x_1\xi_n\xi_{n-1}\cdots\xi_{1}\}\nonumber \\
\G_{n}&:&\{R((n+2)\pi_1),x_1^{n+2}\},\nonumber \\
&&\{R((n+1)\pi_1+\pi_{n-1}),x_1^{n+1}\xi_n\},\{R(n\pi_1),x_1^{n}\Phi\},\cdots,\nonumber \\
&&\{R(4\pi_1),x_1^3\xi_n\cdots\xi_2\},\nonumber \\
&&\{R(2\pi_1+\pi_2),x_1^2\xi_n\cdots\xi_3\Phi\},\{R(2\pi_1),x_1^2\xi_n\xi_{n-1}\cdots\xi_{1}\}\nonumber \\
&\qquad\vdots\nonumber 
\end{eqnarray}

$\Lambda^2(\G_{-1}^*)\cong\{R(2\pi_{n-1}),f_n^2\}\oplus\{R(\pi_1+\pi_{n-1}),f_n\theta_1\}\oplus\{R(\pi_2),\theta_1\theta_2\}\oplus\{R(0),\sum_{i=1}^nf_i\theta_i\}$, where as before $f_i$ and $\theta_i$ are the corresponding dual basis to $x_i$ and $\xi_i$
, respectively.

For $n=3$, there is a $gl_3$-invariant vector in $\Lambda^2(\G_{-1}^*)\otimes\G_{-1}$, namely $\sum_{i=1}^n\theta_i^*\otimes x_i$, where $\theta_i^*$ is the Hodge dual to $\theta_i$.  However, this vector is odd.  For $n\ge 3$, there are no non-zero $gl_n
$-invariant vectors in $\Lambda^2(\G_{-1}^*)\otimes\G_{-1}$ by inspection of the table above.

For $n\ge 3$ there are four linearly independent $sl_n$-invariant cochains in $\Lambda^2(\G_{-1}^*)\otimes\G_0$.  However they are all odd.

Consider now the $sl_n$-invariant cochains in $\Lambda^2(\G_{-1}^*)\otimes\G_2$.  They are all linear combinations of the following vectors:

\begin{eqnarray}
c_1&=&\sum_{i,j}f_if_j\otimes x_i x_j\Phi,\nonumber \\
c_2&=&{1\over 2}\sum_{i=1}^{n-1}(f_i\theta_i-f_{i+1}\theta_{i+1})\otimes(x_i\xi_i-x_{i+1}\xi_{i+1})\Phi \nonumber \\
&+&\sum_{i\not=j}f_i\theta_j\otimes x_i\xi_j\Phi,\nonumber \\
c_3&=&\sum_{i<j}\theta_i\theta_j\otimes\xi_i\xi_j\Phi.\nonumber 
\end{eqnarray}
Compute
\begin{eqnarray}
dc_1(x_1,x_1,x_1)&=&-3[x_1,x_1^2\Phi]=-3x_1^3\not=0,\nonumber \\
dc_2(x_1,x_1,x_1)&=&0,\nonumber \\
dc_3(x_1,x_1,x_1)&=&0,\nonumber \\
dc_1(\xi_1,\xi_2,\xi_3)&=&0,\nonumber \\
dc_2(\xi_1,\xi_2,\xi_3)&=&0,\nonumber \\
dc_3(\xi_1,\xi_2,\xi_3)&\not=&0,\nonumber \\
dc_2(\xi_1,x_1,x_1)&=&[x_1,x_1\xi_1-x_2\xi_2]=x_1\not=0.\nonumber 
\end{eqnarray}
Therefore there are no cocycles in the span of $c_1$, $c_2$ and $c_3$.  Hence there are no $gl_n$-invariant Spencer cocycles in $\Lambda^2(\G_{-1}^*)\otimes\G_2$.

Consider now the $sl_n$-invariant cochains in $\Lambda^2(\G_{-1}^*)\otimes\G_{n-2}$.  They are all linear combinations of the following vectors:
\begin{eqnarray}
c_1 &=& \sum_{i,j}f_if_j\otimes x_i\xi_i^*,\nonumber \\
c_2 &=& (\sum_{i=1}^n f_i\theta_i)\otimes\xi_1\xi_2\cdots\xi_n,\nonumber
\end{eqnarray}
It is evident that $\Phi$ acts non-trivially on $c_1$ and $c_2$, hence these vectors are not $gl_n$-invariant.
Hence there are no $gl_n$-invariant non-trivial Spencer cocycles in $\Lambda^2(\G_{-1}^*)\otimes\G_{n-2}$.

Finally there is an $sl_n$-invariant cocycle in $\Lambda^2(\G_{-1}^*)\otimes\G_n$, namely
$$c=\sum_{i,j}f_if_j\otimes x_ix_j\xi_1\xi_2\cdots\xi_n.$$
But obviously $c$ is not $gl_n$-invariant.

Next consider $HO(4,4)$.  The $sl_4$-invariant cochains in $\Lambda^2(\G_{-1}^*)\otimes\G_0$ and $\Lambda^2(\G_{-1}^*)\otimes\G_4$ are as in the general case and the same argument applies.

There are five cochains now in $\Lambda^2(\G_{-1}^*)\otimes\G_2$, namely
\begin{eqnarray}
c_1&=&\sum_{i,j}f_if_j\otimes x_i x_j\Phi,\nonumber \\
c_2&=&{1\over 2}\sum_{i=1}^{3}(f_i\theta_i-f_{i+1}\theta_{i+1})\otimes(x_i\xi_i-x_{i+1}\xi_{i+1})\Phi \nonumber \\
&+&\sum_{i\not=j}f_i\theta_j\otimes x_i\xi_j\Phi,\nonumber \\
c_3&=&\sum_{i<j}\theta_i\theta_j\otimes\xi_i\xi_j\Phi,\nonumber \\
c_4&=&\sum_{i,j}f_if_j\otimes x_i\xi_j^*,\nonumber \\
c_5&=&(\sum_{i=1}^4 f_i\theta_i)\otimes\xi_1\xi_2\xi_3\xi_4.\nonumber 
\end{eqnarray}

(Recall that $\xi_j^*$ denotes the Hodge dual as usual.) However, $c_4$ and $c_5$ are not $gl_4$-invariant, and we have seen earlier in the general case that there is no cocycle in the span of $c_1$, $c_2$ and $c_3$.

Finally consider $HO(2,2)$.  In this case the $sl_2$-module structure of $HO(2,2)$ is as follows:
\begin{eqnarray}
\G_{-1}&:&\{R(1), x_1\},\{R(1),\xi_2\}\nonumber \\
\G_{0}&:&\{R(2),x_1^2\},\{R(2),x_1\xi_2\},\{R(0),\Phi\},\{R(0),\xi_1\xi_2\}\nonumber \\
\G_{1}&:&\{R(3),x_1^3\},\{R(3),x_1^2\xi_2\},\{R(1),x_1\Phi\},\nonumber \\
&&\{R(1),x_1\xi_1\xi_2\}\nonumber \\
\G_{2}&:&\{R(4),x_1^4\},\{R(4),x_1^3\xi_2\},\{R(2),x_1^2\Phi\},\{R(2),x_1^2\xi_1\xi_2\}\nonumber \\ \qquad\vdots\nonumber \\
\G_{k}&:&\{R(k+2),x_1^{k+2}\},\{R(k+2),x_1^{k+1}\xi_2\},\nonumber \\
&&\{R(k),x_1^{k}\Phi\},\{R(k),x_1^{k}\xi_1\xi_2\}\nonumber \\
&&\qquad\vdots\nonumber 
\end{eqnarray}
$\Lambda^2(\G_{-1}^*)\cong\{ R(2),f_2^2\} \oplus\{ R(2),f_2\theta_1\} \oplus\{ R(0),f_1\theta_1+f_2\theta_2\}\oplus\{ R(0),\theta_1\theta_2\}$.  Here as usual $f_i(x_j)=\delta_{ij}$ and $\theta_i(\xi_j)=\delta_{ij}$, $i,j=1,2$.

In $\Lambda^2(\G_{-1}^*)\otimes\G_0$ the even $sl_2$-invariant cochains are in the span of the following vectors:
\begin{eqnarray}
c_1&=&f_1f_2\otimes(x_1\xi_1-x_2\xi_2)\nonumber \\
&-&f_1^2\otimes x_1\xi_2+f_2^2\otimes x_2\xi_1,\nonumber \\
c_2&=&(f_1\theta_1-f_2\theta_2)\otimes x_1x_2\nonumber \\
&-&f_1\theta_2\otimes x_1^2+f_2\theta_1\otimes x_2^2,\nonumber \\
c_3&=&(f_1\theta_1+f_2\theta_2)\otimes\xi_1\xi_2,\nonumber \\
c_4&=&\theta_1\theta_2\otimes(x_1\xi_1+x_2\xi_2).\nonumber 
\end{eqnarray}

However, none of them is $\Phi$-invariant, as is easily seen.

In $\Lambda^2(\G_{-1}^*)\otimes\G_2$ the even $sl_2$-invariant cochains are in the span of the following vectors:
\begin{eqnarray}
c_1&=&\sum_{i,j}f_if_j\otimes x_i x_j\Phi,\nonumber \\
c_2&=&(f_1\theta_1-f_2\theta_2)\otimes x_1x_2\xi_1\xi_2-f_1\theta_2\otimes
x_1^2\xi_1\xi_2 +f_2\theta_1\otimes x_2^2\xi_1\xi_2.\nonumber 
\end{eqnarray}
Compute
\begin{eqnarray}
dc_1(x_1,x_1,x_1)&=&-3[x_1,x_1^2\Phi]\not=0,\nonumber \\
dc_2(x_1,x_1,x_1)&=&0,\nonumber \\
dc_2(\xi_1,\xi_2,x_1)&=&-3x_1\xi_1\xi_2\not=0.\nonumber
\end{eqnarray}

Hence there are no even $gl_2$-invariant cocycles in $\Lambda^2(\G_{-1}^*)\otimes\G_2$.  Thus we have proved

\proclaim{Proposition 4.5.}
$$H^{l,2}(HO(n)_{-1};HO(n))^{gl_n}_{\bar{0}}=0,\quad l\ge 1;\ n\ge 2.$$

Recall that $HO(n,n)$ is the full prolongation of $\G_{\le 0}$.  Hence by Corollary 2.3 we obtain

\proclaim{Theorem 4.2.} $HO(n,n)$ has no non-trivial filtered deformations for $n\ge 2$.

\remark{Remark 4.1.} Let $\G$ be either $\hat{SHO}(n,n)$, $\hat{SHO}'(n,n)$ or $\hat{HO}(n,n)$ and let $\A$ be its maximal reductive subalgebra with respect to which we have decomposed $\G$.  For $n\not=3$ it is clear from the $\A$-module structure of $\G
$ that the $\A$-module $\G_0^*\otimes\G_{-1}$ has no trivial $\A$-component in its $\A$-module decomposition, and hence it has no trivial $\G_0$-subquotient in its $\G_0$-composition series.  Therefore $H^1(\G_0;\G_{-1})=0$. 
It follows from Proposition 2.7 that $\G$ has no filtered deformation $L$ such that $L_0$ is a maximal subalgebra.  In the case when $n=3$, the vector $1$ and the trivial $sl_n$-module in $\G_0^*\otimes \G_{-1}$ have opposite parity, and so Proposition 2.
7 takes care of this case as well.  In fact it can be shown directly that $\hat{HO}(n,n)$ has no non-trivial filtered deformations at all.  However the remaining two cases do possess non-trivial filtered deformations which turn out to be interesting.  Thi
s is the reason why we have calculated Spencer $2$-cocycles of $\hat{SHO}(n,n)$ and $\hat{SHO}'(n,n)$ in Propositions 4.3 and 4.4.

In the next two remarks we will deal with the Lie superalgebras $SHO(n,n)$\newline $+\C\Phi$, $SHO'(n,n)+\C\Phi$, $\hat{SHO}(n,n)+\C\Phi$ and $\hat{SHO}'(n,n)+\C\Phi$.

\remark{Remark 4.2.} Let $\G$ be either $\hat{SHO}(n,n)+\C\Phi$ or $\hat{SHO}'(n,n)+\C\Phi$.  Here the maximal reductive subalgebra $\A$ of $\G_0$ is of course $gl_n=sl_n+\C\Phi$.  As in Remark 4.1 it follows from Proposition 2.7 that $\G$ has no filtered
 deformation $L$ such that $L_0$ is a maximal subalgebra.  Actually one can show that $\hat{SHO}(n,n)+\C\Phi$ has a unique non-trivial (non-simple) filtered deformation, while $\hat{SHO}'(n,n)+\C\Phi$ has no non-trivial filtered deformations at all.

\remark{Remark 4.3.} Let $\G$ be either $SHO(n,n)+\C\Phi$ or $SHO'(n,n)+\C\Phi$.  Our computations of Spencer $2$-cocycles of $HO(n,n)$, $SHO(n,n)$ and $SHO'(n,n)$ above also show that $\G$ has no non-trivial Spencer $2$-cocycles.  Namely, the computation
 for $HO(n,n)$ shows that all Spencer $2$-cocycles of $\G$ of degree $2$ are odd and so cannot give rise to filtered deformations.  Also it is easy to check that the unique non-trivial Spencer $2$-cocycle of $SHO'(n,n)$ of degree $n$ cannot be invariant u
nder the action of $\Phi$.  Thus if $\{\mu_1,\mu_2,\cdots\}$ is a defining sequence of a filtered deformation $\G_{\epsilon}$, then we may assume that $\mu_i$ vanishes for all $i$ when restricted to $\G_{-1}$.  In particular it follows that $[\G_{-1},\G_{
-1}]_{\epsilon}$ has trivial projection onto $\C\Phi$.  Now for $n\not=3$ the trivial $sl_n$-module does not appear in $\G_0\otimes\G_{-1}$ and hence $[\G_0,\G_{-1}]_{\epsilon}$ projects trivially onto $\C\Phi$.  When $n=3$, $[\G_0,\G_{-1}]_{\epsilon}$ pr
ojects trivially onto $\C\Phi$ due to parity reason.  Of course $[\G_1,\G_{-1}]_{\epsilon}$ projects trivially onto $\C\Phi$.  Therefore $[\G,\G]_{\epsilon}$ projects trivially onto $\C\Phi$ and hence $\G_{\epsilon}$ cannot be simple.  Again here one can 
determine all non-trivial filtered deformations.  It turns out that $SHO(n,n)+\C\Phi$ has a unique non-trivial (non-simple) filtered deformation, while $SHO'(n,n)+\C\Phi$ has no non-trivial filtered deformation at all.

We shall next consider $SKO'(n,n+1;\beta)$ and $SKO(n,n+1;\beta)$.  Here the even part of $\G_0$ is $gl_n$.  As an $sl_n$-module $SKO'(n,n+1;\beta)$, $\beta\not=1$, decomposes as follows (as usual we include a highest weight vector and $\Phi=\sum_{i=1}^nx
_i\xi_i$):
\begin{eqnarray}
\G_{-2}&:&\{R(0),1\}\G_{-1}:\{R(\pi_1), x_1\},\{R(\pi_{n-1}),\xi_n\}\nonumber \\
\G_{0}&:&\{R(2\pi_1),x_1^2\},\{R(\pi_1+\pi_{n-1}),x_1\xi_n\},\nonumber \\
&&\{R(0),\tau+\beta\Phi\},\{R(\pi_{n-2}),\xi_n\xi_{n-1}\}\nonumber \\
\G_{1}&:&\{R(3\pi_1),x_1^3\},\{R(2\pi_1+\pi_{n-1}),x_1^2\xi_n\},\nonumber \\
&&\{R(\pi_1),x_1(\tau+{{\beta n-1}\over{n+1}}\Phi)\},\nonumber \\
&&\{R(\pi_1+\pi_{n-2}),x_1\xi_{n}\xi_{n-1}\},\nonumber \\
&&\{R(\pi_{n-1}),\xi_n(\tau+{{\beta n-1}\over{n-1}}\Phi)\},\nonumber \\
&&\{R(\pi_{n-3}),\xi_n\xi_{n-1}\xi_{n-2}\}\nonumber \\
\G_{2}&:&\{R(4\pi_1),x_1^4\},\{R(3\pi_1+\pi_{n-1}),x_1^3\xi_n\},\nonumber \\
&&\{R(2\pi_1),x_1^2(\tau+{{\beta n-2}\over{n+2}}\Phi)\},\nonumber \\
&&\{R(2\pi_1+\pi_{n-2}),x_1^2\xi_{n}\xi_{n-1}\},\nonumber \\
&&\{R(\pi_1+\pi_{n-1}),x_1\xi_n(\tau+{{\beta n-2}\over{n}}\Phi)\},\nonumber \\
&&\{R(\pi_1+\pi_{n-3}),x_1\xi_n\xi_{n-1}\xi_{n-2}\},\nonumber \\
&&\{R(\pi_{n-2}),\xi_{n-1}\xi_{n-2}(\tau+{{\beta n-2}\over{n-2}}\Phi)\},\nonumber \\
&&\{R(\pi_{n-4}),\xi_n\xi_{n-1}\xi_{n-2}\xi_{n-3}\},\nonumber \\
&&\qquad\vdots,\nonumber \\
\G_{n-2}&:&\{R(n\pi_1),x_1^n\},\{R((n-1)\pi_1+\pi_{n-1}),x_1^{n-1}\xi_n\},\nonumber \\
&&\{R((n-2)\pi_1),x_1^{n-2}(\tau+{{\beta n-n+2}\over{2n-2}}\Phi)\}\cdots,\nonumber \\
&&\{R(2\pi_1),x_1\xi_n\cdots\xi_2\},\nonumber \\
&&\{R(\pi_2),\xi_n\cdots\xi_3(\tau+{{\beta n-n+2}\over{2}}\Phi)\},\nonumber \\
&&\{R(0),\xi_n\xi_{n-1}\cdots\xi_{1}\}\nonumber \\
\G_{n-1}&:&\{R((n+1)\pi_1),x_1^{n+1}\},\{R(n\pi_1+\pi_{n-1}),x_1^n\xi_n\},\nonumber \\
&&\{R((n-1)\pi_1),x_1^{n-1}(\tau+{{\beta n-n+1}\over{2n-1}}\Phi)\},\cdots,\nonumber \\ &&\{R(3\pi_1),x_1^2\xi_n\cdots\xi_2\},\nonumber \\
&&\{R(\pi_1+\pi_2),x_1\xi_n\cdots\xi_3(\tau+{{\beta n-n+1}\over{3}}\Phi)\},\nonumber \\
&&\{R(\pi_1),\xi_n\xi_{n-1}\cdots\xi_{2}(\tau+{{\beta n-n+2}\over{1}}\Phi)\}\nonumber \\
\G_{n}&:&\{R((n+2)\pi_1),x_1^{n+2}\},\{R((n+1)\pi_1+\pi_{n-1}),x_1^{n+1}\xi_n\},\nonumber \\
&&\{R(n\pi_1),x_1^{n}(\tau+{{\beta n-n}\over{2n}}\Phi)\},\cdots,\nonumber \\
&&\{R(4\pi_1),x_1^3\xi_n\cdots\xi_2\},\{R(2\pi_1+\pi_2),x_1^2\xi_n\cdots\xi_3(\tau+{{\beta n-n}\over{4}}\Phi)\},\nonumber \\
&&\{R(2\pi_1),x_1\xi_n\xi_{n-1}\cdots\xi_{2}(\tau+{{\beta n-n}\over{2}}\Phi)\}\nonumber \\
\G_{n+1}&:&\{R((n+3)\pi_1),x_1^{n+3}\},\{R((n+2)\pi_1+\pi_{n-1}),x_1^{n+2}\xi_n\},\nonumber \\
&&\{R((n+1)\pi_1),x_1^{n+1}(\tau+{{\beta n-n-1}\over{2n+1}}\Phi)\},\cdots,\{R(5\pi_1),x_1^4\xi_n\cdots\xi_2\},\nonumber \\
&&\{R(3\pi_1+\pi_2),x_1^3\xi_n\cdots\xi_3(\tau+{{\beta n-n-1}\over{5}}\Phi)\},\nonumber \\
&&\{R(3\pi_1),x_1^2\xi_n\xi_{n-1}\cdots\xi_{2}(\tau+{{\beta n-n-1}\over{3}}\Phi)\}\nonumber \\
&&\qquad\vdots\nonumber 
\end{eqnarray}
In the case when $\beta=1$ an extra component $\{R(0),\tau\xi_n\xi_{n-1}\cdots\xi_{1}\}$ is included in $\G_n$.  The structure of $SKO(n,n+1;\beta)$ is then easily derived from (3.2).
$\Lambda^2(\G_{-}^*)\cong\{R(2\pi_{n-1}),f_n^2\}\oplus\{R(\pi_1+\pi_{n-1}),f_n\theta_1\}\oplus\{R(\pi_2),\theta_1\theta_2\}\oplus\{R(\pi_{n-1}),1^*f_n\}\oplus\{R(\pi_1),1^*\theta_1\}\oplus\{R(0),1^*1^*\}\oplus\{R(0),\sum_{i=1}^nf_i\theta_i\}$, where $f_i$
, $\theta_i$ and $1^*$ are the corresponding dual basis to $x_i$, $\xi_i$ and $1$, respectively.  Below we will find $gl_n$-invariant $2$-cocycles.  Just as for $HO(n,n)$ the cases $n=2,4$ again need to be considered separately, however, the analysis is c
ompletely analogous and we will omit these cases.  For our calculations below we shall need the following lemma, whose proof is straightforward.

\proclaim{Lemma 4.1.} Let $f$ be a monomial in $\C[x_1,\cdots,x_n,\xi_1,\cdots,\xi_n]$ and let $\lambda\in\C$.  Let $o(f)$ and $e(f)$ denote the number of $\xi_i$'s and $x_i$'s in $f$, respectively. Then

\indent(i) $[\tau+\beta\Phi,(\tau+\lambda\Phi)f]=((1-\beta)o(f)+(1+\beta)e(f))(\tau+\lambda\Phi)f$.

\indent(ii) $[(\tau+\beta\Phi),f]=(-2+(1-\beta)o(f)+(1+\beta)e(f))f$.

For $n=3$ there is an $sl_n$-invariant $2$-cochain of degree $1$.  However, as in the case of $HO(3,3)$, this vector is odd.

The $sl_n$-invariant Spencer $2$-cochains of degree $2$ are all odd, and so are excluded.

There are six linearly independent $sl_n$-invariant Spencer $2$-cochains of degree $4$.  They are as follows:
\begin{eqnarray}
c_1&=&\sum_{i,j}f_if_j\otimes x_ix_j(\tau+{{\beta n-2}\over{n+2}}\Phi),\nonumber \\
c_2&=&{1\over 2}\sum_{i=1}^n(f_i\theta_i-f_{i+1}\theta_{i+1})\otimes(x_i\xi_i-x_{i+1}\xi_{i+1})(\tau+{{\beta n-2}\over{n}}\Phi)\nonumber \\
&+&\sum_{i\not=j}f_i\theta_j\otimes x_i\xi_j(\tau+{{\beta n-2}\over{n}}\Phi),\nonumber \\
c_3&=&\sum_{i<j}\theta_i\theta_j\otimes\xi_i\xi_j(\tau+{{\beta n-2}\over{n-2}}\Phi),\nonumber \\
c_4&=&\sum_{i=1}^n1^*f_i\otimes x_i(\tau+{{\beta n-1}\over{n+1}}\Phi),\nonumber \\
c_5&=&\sum_{i=1}^n1^*\theta_i\otimes\xi_i(\tau+{{\beta n-1}\over{n-1}}\Phi),\nonumber \\
c_6&=&1^*1^*\otimes(\tau+\beta\Phi).\nonumber 
\end{eqnarray}
However, by Lemma 4.1 all these vectors have $(\tau+\beta\Phi)$-eigenvalue $4$, and hence are not $gl_n$-invariant.

We have three linearly independent $sl_n$-invariant vectors of degree $n$, namely
\begin{eqnarray}
c_1&=&\sum_{i,j}f_if_j\otimes x_i\xi_j^*,\nonumber \\
c_2&=&\sum_{i=1}^nx_i\theta_i\otimes\xi_1\xi_2\cdots\xi_n,\nonumber \\
c_3&=&\sum_{i=1}^n1^*f_i\otimes\xi_i^*,\nonumber 
\end{eqnarray}
\noindent where as usual $\xi_i^*$ is the Hodge dual of $\xi_i$. By Lemma 4.1 $\tau+\beta\Phi$ acts on these vectors as the scalar $(1-\beta)n$ and hence they are $gl_n$-invariant if and only if $\beta=1$.  However there are two $1$-cochains of degree $n$
, namely
\begin{eqnarray}
b_1&=&\sum_{i=1}^nf_i\otimes\xi_i^*(\tau+{{\beta n-n+1}\over{1}}\Phi),\nonumber \\
b_2&=&1^*\otimes \xi_1\xi_2\cdots\xi_n,\nonumber 
\end{eqnarray}
which are $gl_n$-invariant if and only if $\beta=1$.  It is easy to check that $db_1$ and $db_2$ are linearly independent and lie in the span of $c_1$, $c_2$ and $c_3$.  But $dc_2(\xi_1,x_1,x_1)=-2\xi_i^*\not=0$ and thus there are no non-trivial cocycles 
in the span of $c_1$, $c_2$ and $c_3$.

Of degree $n+2$ there are three $sl_n$-invariant cochains, namely
\begin{eqnarray}
c_1&=&\sum_{i,j}f_if_j\otimes x_i\xi_j^*(\tau+{{\beta n-n}\over{2}}\Phi),\nonumber \\
c_2&=&\sum_{i=1}^n1^*f_i\otimes\xi_i^*(\tau+{{\beta n-n+1}\over{1}}\Phi),\nonumber \\
c_3&=&1^*1^*\otimes\xi_1\xi_2\cdots\xi_n.\nonumber 
\end{eqnarray}
Note that they are even cochains if and only if $n$ is odd.
By Lemma 4.1 they are $\tau+\beta\Phi$-invariant if and only if $\beta={{n+2}\over n}$.
Compute ($\beta={{n+2}\over n}$)
\begin{eqnarray}
dc_1(\xi_1,x_1,x_1)&=&\tau\xi_i^*+3x_1\xi_1\xi_2\cdots\xi_n,\nonumber \\
dc_3(\xi_1,x_1,x_1)&=&0,\nonumber \\
dc_1(\xi_1,x_1,1)&=&0,\nonumber \\
dc_3(\xi_1,x_1,1)&=&\xi_1\xi_2\cdots\xi_n.\nonumber
\end{eqnarray}
Hence the space of cocycles in the span of $c_1$, $c_2$ and $c_3$ is at most one.  We will construct a cocycle of degree $n+2$ in the next section.

For $\beta=1$ there exists an $sl_n$-invariant cochain of degree $n+4$, namely
$$c=1^*1^*\otimes\tau\xi_1\xi_2\cdots\xi_n.$$
By Lemma 4.1 its $(\tau+\beta\Phi)$-eigenvalue is $4$, hence it is not $gl_n$-invariant.

We summarize our calculations above.

\proclaim{Proposition 4.6} For $n\ge 2$
\begin{eqnarray}
&&H^{l,2}(SKO(n,n+1;\beta)_{-};SKO(n,n+1;\beta))_{\bar{0}}^{gl_n}=0,\;\forall l{\rm\ and\ }\beta\not={{n+2}\over n},\nonumber \\
&&H^{l,2}(SKO(n,n+1;{{n+2}\over n})_-;SKO(n,n+1;{{n+2}\over n}))_{\bar{0}}^{gl_n}=0,\; l\not={n+2},\nonumber \\
&&H^{n+2,2}(SKO(n,n+1;{{{n+2}\over n}})_{-};SKO(n,n+1;{{{n+2}\over n}}))_{\bar{0}}^{gl_n}=\C, n{\rm\ odd},\nonumber \\
&&H^{n+2,2}(SKO(n,n+1;{{{n+2}\over n}})_{-};SKO(n,n+1;{{{n+2}\over n}}))_{\bar{0}}^{gl_n}=0, n{\rm\ even}.\nonumber 
\end{eqnarray}

\proclaim{Theorem 4.3.} The Lie superalgebras $SKO(n,n+1;\beta)$ and $SKO'(n,n+1;\beta)$, for $\beta\not={{n+2}\over n}$ or $n$ even, have no non-trivial filtered deformations.

\proof We shall always assume that $\beta\not={{n+2}\over n}$.  We know that  $SKO'(n,n+1;\beta)$ is a full prolongation, and hence by Proposition 4.6 and Corollary 2.3, it has no non-trivial filtered deformations.  Since $SKO(n,n+1;\beta)=SKO'(n,n+1;\beta)$ for $\beta\not=1,{{n-2}\over n}$, we are left to consider two cases.

Now $SKO(n,n+1;{{n-2}\over n})_{(1)}$ contains no trivial $sl_n$-module.  Thus it is an almost full prolongation by (3.2) and hence Proposition 4.6 and Corollary 2.3 take care of this case as well.

Now consider $\G=SKO(n,n+1;1)$.  In this case it is not an almost full prolongation.  We need to go back to the proof of Proposition 2.6, from which and (3.2) it follows that if $L$ is a filtered deformation of $\G$, then $L$ can be given a defining seque
nce $\{\mu_1,\mu_2,\cdots\}$ with the properties that $\mu_i|_{\G_-\times\G}=0$ for $i<n$, $\mu_{n}(\G_-,a)=0$ for $a\in\G_0$ not lying in the trivial $sl_n$-component and $\mu_{n}(x,\tau+\Phi)=\lambda[\tau\xi_1\xi_2\cdots\xi_n,x]$, for all $x\in\G_-$ and
 for some $\lambda\in\C$.  Furthermore by Proposition 4.6 we may also assume that $\mu_n$ vanishes when restricted to $\G_-$. Of course this only makes sense if $n$ is even.  Hence we may assume that $n$ is an even integer from now on. For a fixed $i$ we 
let $b=\xi_i(\tau+\Phi)\in\G_1$, lying in the irreducible $sl_n$-module $R(\pi_{n-1})$.  We then have for $x\in\G_-$:
$$[b,x]_{\epsilon}=[b,x]+\mu_{n}(b,x)\epsilon^{n}+\cdots.$$
Taking $x\in\G_{-1}$ we have $\mu_{n}(b,x)\subset\G_{n}$.  But the irreducible $sl_n$-modules $R(0)$, $R(\pi_1+\pi_{n-1})$ and $R(\pi_{n-2})$ do not appear in $\G_{n}$.  Hence $\mu_{n}(b,\G_{-1})=0$ by Proposition 2.4.  Using this fact we compute the Jaco
bi identity of the triple $x_i,\xi_j,\xi_i(\tau+\Phi)$ for $i\not=j$ and derive that $\lambda[\tau\xi_1\xi_2\cdots\xi_n,x_j]=0$.  Thus $\lambda=0$.  Therefore $\mu_{n}|_{\G_-\times\G_0}=0$.  Now $\G_{(1)}$ only has one more trivial $sl_n$-component, namel
y $\xi_1\xi_2\cdots\xi_n$.  However it has parity different from that of $\tau\xi_1\xi_2\cdots\xi_n$.  Thus we conclude that $\mu_{n}|_{\G_-\times\G}=0$.  Now $\G$ is a full prolongation of degree $n+1$, which combined with Proposition 4.6, allows us to a
pply Proposition 2.6 to conclude that $L$ is a trivial filtered deformation. {$\qed$\medskip}

We shall now consider the Lie superalgebra $\G=H(2n,m)$.  Here $(\G_0)_{\bar{0}}$ is isomorphic to $\A=sp_{2n}\oplus so_m$.  With respect to $\A$, $\G$ decomposes as follows ($\pi_i$ and $\tilde{\pi_i}$ are the respective fundamental weights of $sp_{2n}$ 
and $so_m$.):

\begin{eqnarray}
\G_{-1}:\quad&&R(\pi_1), R(\tilde{\pi}_1)\nonumber \\
\G_{0}:\quad&&R(2\pi_1), R(\pi_1)\otimes R(\tilde{\pi}_1), R(\tilde{\pi}_2)\nonumber \\
&&\quad\vdots\qquad\vdots\nonumber \\
\G_{k}:\quad&& R((k+2)\pi_1), R((k+1)\pi_1)\otimes R(\tilde{\pi}_1),\cdots, R(\tilde{\pi}_{k+2})\nonumber \\
\G_{m-2}:\quad&&R(m\pi_1), R((m-1)\pi_1)\otimes R(\tilde{\pi}_1),\cdots,R(\pi_1)\otimes R(\tilde{\pi}_{m-1}), R(0)\nonumber \\
&&\quad\vdots\qquad\vdots\nonumber \\
\G_{s}:\quad&&R((s+2)\pi_1), R((s+1)\pi_1)\otimes R(\tilde{\pi}_1),\cdots,
\nonumber \\ && R((s-m+3)\pi_1)\otimes R(\tilde{\pi}_{m-1}), R((s-m+2)\pi_1)\nonumber \\
&&\quad\vdots\qquad\vdots \nonumber 
\end{eqnarray}

Continuing using the notation adapted in Section 3, we let $dp_i$, $dq_i$ for $i=1,\cdots,n$ and $d\xi_j$ for $j=1,\cdots,m$ denote the dual basis to $p_i$, $q_i$ and $\xi_j$, respectively.  Then as $sp_{2n}\oplus so_m$-modules we have 
$$\Lambda^2(\G_{-1}^*)\cong R(\pi_2)\oplus R(0)\oplus R(\pi_1)\otimes R(\tilde{\pi}_1)\oplus R(2\tilde{\pi}_1)\oplus R(0),\eqno{(4.1)}$$
where the two trivial components are spanned by the vectors $\sum_{i=1}^n dp_idq_i$ and $\sum_{i=1}^m d\xi_i^2$, respectively.  We need to make some further clarifications of (4.1):  $R(\pi_2)=\emptyset$ if $n=1$.  Also $R(2\tilde{\pi}_1)$ is understood a
s follows: It is irreducible of highest weight $2\tilde{\pi}_1$ only when $m\ge 5$.  For $m=4$ it is isomorphic to $R(2)\otimes \hat{R}(2)$, where $so_4\cong sl_2\oplus\hat{sl}_2$.  For $m=3$ it is $R(4)$, where we identify $so_3$ with $sl_2$.  For $m=2$ 
it is isomorphic to $\C_+\oplus \C_-$, where $\xi_1\xi_2$ acts on the one-dimensional spaces $\C_+$ and $\C_-$ as the scalars $2\sqrt{-1}$ and $-2\sqrt{-1}$, respectively.  For $m=1$, it is empty.

Note that all the modules are self-contragredient.  Furthermore $R(\pi_2)$ doesn't appear in $\G_k$ for any $k$.  Also $R(\tilde{\pi}_k)$ are all non-isomorphic except for $R(\tilde{\pi}_{m-i})\cong R(\tilde{\pi}_{i})$.  Finally the component $R(2\tilde{\
pi}_1)$ in (4.1) is not isomorphic to $R(\tilde{\pi}_k)$ in $\G_k$ for any $k$.  From this it follows that the only $sp_{2n}\oplus so_m$-invariant cochains are:

In $\Lambda^2(\G_{-1}^*)\otimes\G_0$:
$$c_0=\sum_{i,j}dp_id\xi_j\otimes p_i\xi_j+dq_id\xi_j\otimes q_i\xi_j.$$
In $\Lambda^2(\G_{-1}^*)\otimes\G_{m-2}$:
\begin{eqnarray}
c_1&=&(\sum_{i=1}^n dp_idq_i)\otimes\xi_1\cdots\xi_m,\nonumber \\
c_2&=&(\sum_{i=1}^m d\xi_i^2)\otimes\xi_1\cdots\xi_m,\nonumber \\
c_3&=&\sum_{i=1}^n dp_id\xi_j\otimes p_i\xi_j^*+dq_id\xi_j\otimes q_i\xi_j^*,\nonumber 
\end{eqnarray}
\noindent where $\xi_j^*$ as usual stands for the Hodge dual of $\xi_j$.

Suppose that $m\not=2$. Now $dc_0(p_1,\xi_1,\xi_1)=2p_1\not=0$.  Hence $c_0$ is not a cocycle.  On the other hand
\begin{eqnarray}
dc_1(\xi_1,\xi_1,\xi_1)&=&0,\nonumber \\
dc_2(\xi_1,\xi_1,\xi_1)&=&-3\xi_1^*\not=0,\nonumber \\
dc_3(\xi_1,\xi_1,\xi_1)&=&0,\nonumber \\
dc_1(p_1,q_1,\xi_1)&=&-\xi_1^*\not=0.\nonumber 
\end{eqnarray}

From this it follows that the space of cocycles in $\Lambda^2(\G_{-1}^*)\otimes\G_{m-2}$ is at most 1-dimensional.  Let
$$b=\sum_{i=1}^n dp_i\otimes p_i\xi_1\cdots\xi_m+dq_i\otimes q_i\xi_1\cdots\xi_m.$$  It is an $sp_{2n}\oplus so_m$-invariant $1$-cochain such that $db\in\Lambda^2(\G_{-1}^*)\otimes\G_{m-2}$ and $db\not=0$.  Hence any $sp_{2n}\oplus so_m$-invariant cocycle
 in $\Lambda^2(\G_{-1}^*)\otimes\G_{m-2}$ must be a coboundary.

Now if $m=2$, then all four cochains appear in $\Lambda^2(\G_{-1}^*)\otimes\G_0$.  We have
$$dc_0(\xi_1,\xi_1,\xi_1)=0.$$
Hence every cocycle must be a linear combination of $c_0$, $c_1$ and $c_3$.  However,
\begin{eqnarray}
dc_0(p_1,\xi_1,\xi_1)&\not=&0,\nonumber \\
dc_1(p_1,\xi_1,\xi_1)&=&0,\nonumber \\
dc_3(p_1,\xi_1,\xi_1)&=&0.\nonumber 
\end{eqnarray}
Thus any cocycle must be a linear combination of $c_1$ and $c_3$.  But we have seen from the general case that it must be a coboundary.

\proclaim{Proposition 4.7.} For $n\ge 1$ and $m\ge 0$
$$H^{l,2}(H(2n,m)_{-1};H(2n,m))^{sp_{2n}\oplus so_m}_{\bar{0}}=0,\quad\forall l.$$

By Corollary 2.3, since $H(2n,m)$ is a full prolongation, we arrive at

\proclaim{Theorem 4.4.} $H(2n,m)$ has no non-trivial filtered deformations for $n\ge 1$ and $m\ge 0$.

\heading{5}{Existence and uniqueness of filtered deformations of \\ $SHO'(n,n)$, $\hat{SHO}(n,n)$ and $SKO(n,n+1;{{n+2}\over {n}})$}

In this section we will construct non-trivial filtered deformations of the Lie superalgebras $\hat{SHO}(n,n)$, $SHO'(n,n)$ and $\hat{SHO}'(n,n)$, for $n$ even, and for $SKO(n,n+1;{{n+2}\over n})$, for $n$ odd.  From Propositions 4.2, 4.3, 4.4 and 4.6 and 
Corollary 2.5 it follows that such filtered deformations are necessarily unique.

Let $\G$ be either $\hat{SHO}(n,n)$ or $\hat{SHO}'(n,n)$.  As usual we identify $\G$ with a subalgebra in $\Lambda(n,n)$ with the odd Poisson bracket.  Recall that $\Lambda(n)$ is naturally $\Z$-graded so that we may write $\Lambda(n)=\oplus_{j=0}^n\Lambda(n)_j$.  We let $\G^{j}=\G\cap(\C[x_1,x_2,\cdots,x_n]\otimes\Lambda(n)_j)$.

\proclaim{Lemma 5.1.} Let $f\in\G^{j}$ with $j\ge 1$, where $\G$ is either $\hat{SHO}(n,n)$, $\hat{SHO}'(n,n)$ or ${SHO}'(n,n)$.  Then $[\xi_1\xi_2\cdots\xi_n,f]=0$.

\proof Note that we have $[\G^{i},\G^{j}]\subset\G^{i+j-1}$.  Hence $[\xi_1\xi_2\cdots\xi_n,\G^{j}]\subset\G^{n+j-1}$.  In particular if $j\ge 1$, $[\xi_1\xi_2\cdots\xi_n,\G^{j}]\subset\G^{n}$.  But $\G^{n}=\C\xi_1\xi_2\cdots\xi_n$ or $\G^{n}=0$.  On the 
other hand we know that $\C\xi_1\xi_2\cdots\xi_n$ doesn't lie in the derived algebra.  Thus $[\xi_1\xi_2\cdots\xi_n,f]=0$. {$\qed$\medskip}

We define a super-skewsymmetric bilinear map $\mu_n:\G\wedge\G\rightarrow\G$ of degree $n$ as follows:
\begin{eqnarray}
\mu_n(f,g)&=&[\xi_1\xi_2\cdots\xi_n,fg]\quad {\rm if}\ f,g\in\G^{0},\nonumber \\
\mu_n(f,g)&=&0\qquad\qquad\qquad\quad{\rm otherwise.}\nonumber \end{eqnarray}

In $\G=\hat{SHO}(n,n)$ or $\G=\hat{SHO}'(n,n)$ we define a new bracket $[\cdot,\cdot]_{\epsilon}$ using this $\mu_n$, i.e. we set
$$[f,g]_{\epsilon}=[f,g]+\mu_n(f,g)\epsilon^n.\eqno{(5.1)}$$

\proclaim{Proposition 5.1.} The deformed bracket (5.1) defines a unique non-trivial filtered deformation of $\hat{SHO}(n,n)$ and $\hat{SHO}'(n,n)$.

\proof It suffices to check the $[\cdot,\cdot]_{\epsilon}$ is a Lie bracket.  Since $[\cdot,\cdot]_{\epsilon}$ is obviously super-skewsymmetric, we only need to verify that the Jacobi identity is satisfied.

Let $\G=\hat{SHO}(n,n)$ or $\G=\hat{SHO}'(n,n)$. Recall that
$[\cdot,\cdot]$ is a Poisson bracket, which means that we have
$$[h,fg]=[h,f]g+(-1)^{p(h)(p(f)+1)}f[h,g].\eqno{(5.2)}$$
Here and
below we mean by $p(f)$ the parity of $f$ as an element of the Lie
superalgebra.  So for example $p(x_i)=1$ and $p(\xi_i)=0$.  We remind
the reader again that the Lie superalgebra $\hat{HO}(n,n)$ is
isomorphic to the odd Poisson superalgebra $\Lambda(n,n)$ with
reversed parity, which accounts for the additional $+1$ in (5.2).  We
need to check
$$[h,[f,g]_{\epsilon}]_{\epsilon}=[[h,f]_{\epsilon},g]_{\epsilon}+p(h,f)[f,[h,g]_{\epsilon}]_{\epsilon}.\eqno{(5.3)}$$
By Lemma 5.1 $\mu_n(a,\mu_n(b,c))=0$, for $a,b,c\in\G$.  Thus the
left-hand side of (5.3) is
$$[h,[f,g]]+(\mu_n(h,[f,g])+[h,\mu_n(f,g)])\epsilon^n,$$ while the
right hand side of (5.3) is
\begin{eqnarray}
[[h,f],g]&+&p(h,f)[f,[h,g]]\nonumber \\
&+&(\mu_n([h,f],g)+[\mu_n(h,f),g]\nonumber \\
&+&p(h,f)\mu_n(f,[h,g])+p(h,f)[f,\mu_n(h,g)])\epsilon^n.\nonumber 
\end{eqnarray}

Thus (5.3) is equivalent to
\begin{eqnarray}
 &&mu_n(h,[f,g])+[h,\mu_n(f,g)]={(5.4)}\nonumber \\
&&\mu_n([h,f],g)+[\mu_n(h,f),g]+p(h,f)\mu_n(f,[h,g])+p(h,f)[f,\mu_n(h,g)].\nonumber
\end{eqnarray}

So we need to verify (5.4) for $f,g,h\in\G$.  It is easy to see that
if one of the $f$, $g$ or $h$ lies in $\G^{j}$, for $j\ge 2$, then the
left-hand side and the right-hand side of (5.4) are zero.  Thus we may
assume that $f,g,h\in\G^{j}$, $j=0,1$.  Now it is as easy to see that
if any two of the $f,g,h$ lie in $\G^{1}$, then (5.4) gives again
$0=0$.  Hence we may assume that either $f,g,h\in\G^{0}$ or exactly
one of them lies in $\G^{1}$ and while the other two lie in $\G^{0}$.
We will consider those cases separately.

\noindent{\csc Case 1.} $f,g,h\in\G^{0}$.  In this case, noting that $p(f)=p(g)=p(h)=1$ and $[\G^{0},\G^{0}]=0$ (5.4) reduces to
$$[h,[\xi_1\cdots\xi_n,fg]]=[[\xi_1\cdots\xi_n,hf],g]-[f,[\xi_1\cdots\xi_n,hg]].\eqno{(5.5)}$$
Now the left-hand side of (5.5) equals
\begin{eqnarray}
&&[h,[\xi_1\cdots\xi_n,f]g]+[h,f[\xi_1\cdots\xi_n,g]]=\nonumber \\
&&[h,[\xi_1\cdots\xi_n,f]]g+f[h,[\xi_1\cdots\xi_n,g]].\nonumber 
\end{eqnarray}
The right-hand side of (5.5) is
\begin{eqnarray}
&&[[\xi_1\cdots\xi_n,h]f,g]+[h[\xi_1\cdots\xi_n,f],g]-[f,[\xi_1\cdots\xi_n,h]g]\nonumber \\&&-[f,h[\xi_1\cdots\xi_n,g]] =[[\xi_1\cdots\xi_n,h],g]f+h[[\xi_1\cdots\xi_n,f],g]\nonumber \\ &&-[f,[\xi_1\cdots\xi_n,h]]g- h[f,[\xi_1\cdots\xi_n,g]] . \nonumber 
\end{eqnarray}
Thus we are to show
\begin{eqnarray}
[h,[\xi_1\cdots\xi_n,f]]g+f[h,[\xi_1\cdots\xi_n,g]]&=&
%{\;\qquad \qquad \qquad \qquad \qquad\qquad(5.6)}\nonumber \\
[[\xi_1\cdots\xi_n,h],g]f+h[[\xi_1\cdots\xi_n,f],g]\nonumber \\ &-&[f,[\xi_1\cdots\xi_n,h]]g-h[f,[\xi_1\cdots\xi_n,g]].\nonumber \\&&{\;\;\;\;\qquad \qquad \qquad \qquad \qquad\qquad(5.6)}\nonumber
\end{eqnarray}
But (5.6) is equivalent to saying that
$$-[\xi_1\cdots\xi_n,[f,h]]g-[\xi_1\cdots\xi_n,[h,g]]f=h[\xi_1\cdots\xi_n,[f,g]].\eqno{(5.7)}$$But (5.7) is obviously true, since $[\G^{0},\G^{0}]=0$.

\noindent{\csc Case 2.} One of the $f,g,h$ is in $\G^{1}$, while the other two are in $\G^{0}$.  We will assume that $h\in\G^{1}$.  Other cases are analogous.  In this case (5.4) reduces to
$$[h,[\xi_1\cdots\xi_n,fg]]=[\xi_1\cdots\xi_n,[h,f]g]+[\xi_1\cdots\xi_n,f[h,g]].\eqno{(5.8)}$$
(5.8) is equivalent to
\begin{eqnarray}
[h,[\xi_1\cdots\xi_n,f]g]+[h,f[\xi_1\cdots\xi_n,g]]&=&[\xi_1\cdots\xi_n,[h,f]]g+[h,f][\xi_1\cdots\xi_n,g]\nonumber \\ &+&[\xi_1\cdots\xi_n,f][h,g]\nonumber \\ &+&f[\xi_1\cdots\xi_n,[h,g]],\nonumber \end{eqnarray}
which in term is equivalent to
\begin{eqnarray}
&[h,[\xi_1\cdots\xi_n,f]]g+[\xi_1\cdots\xi_n,f][h,g]
+[h,f][\xi_1\cdots\xi_n,g]+f[h,[\xi_1\cdots\xi_n,g]]\nonumber \\
=&[\xi_1\cdots\xi_n,[h,f]]g+[h,f][\xi_1\cdots\xi_n,g]
+[\xi_1\cdots\xi_n,f][h,g]+f[\xi_1\cdots\xi_n,[h,g]]. \nonumber \end{eqnarray}
But this is equivalent to saying that
$$f[[\xi_1\cdots\xi_n,h],g]+[[\xi_1\cdots\xi_n,h],f]g=0,$$
which is certainly true, since $[\xi_1\cdots\xi_n,h]=0$ by Lemma 5.1. {$\qed$\medskip}

Denote these filtered deformations of $\hat{SHO}(n,n)$ and $\hat{SHO}'(n,n)$ by $\hat{SHO}(n,n)_{\epsilon}$ and $\hat{SHO}'(n,n)_{\epsilon}$, respectively.  $\hat{SHO}(n,n)_{\epsilon}$ is a simple Lie superalgebra.  Now in $\hat{SHO}'(n,n)_{\epsilon}$ $1$
 is no longer central.  But it is easy to see that $1-\xi_1\cdots\xi_n$ is central.  Dividing by the ideal $\C(1-\xi_1\cdots\xi_n)$ we obtain a filtered deformation of $SHO'(n,n)$, which we denote by $SHO'(n,n)_{\epsilon}$.  The Lie bracket in $SHO'(n,n)_
{\epsilon}$ is given by:
\begin{eqnarray}
&&[f,g]_{\epsilon} = [\xi_1\cdots\xi_n,fg],\quad f,g\in\G^{0},\nonumber \\
&&[x_i,\xi_j]_{\epsilon} = \delta_{ij}\xi_1\cdots\xi_n,\qquad\qquad\qquad\qquad\qquad\qquad\qquad\qquad\qquad\qquad{(5.9)}\nonumber \\&&[f,g]_{\epsilon} = [f,g],\quad{\rm otherwise}.\nonumber 
\end{eqnarray}

Define a map $\rho:\hat{SHO}(n,n)_{\epsilon}\rightarrow SHO'(n,n)_{\epsilon}$
by
\begin{eqnarray}
\rho(f)&=&f,\qquad f\in\G_{(-1)},\nonumber \\
\rho(1)&=&\xi_1\xi_2\cdots\xi_n.\nonumber 
\end{eqnarray}
It is easy to see that $\rho$ is an isomorphism of Lie superalgebras.  We summarize our discussion above in the following theorem.

\proclaim{Theorem 5.1.} Let $n\ge 2$ be a positive even integer.

\indent (i) $SHO'(n,n)_{\epsilon}$ with bracket given as in (5.9) defines a unique non-trivial filtered deformation of $SHO'(n,n)$.  It is a simple filtered Lie superalgebra.

\indent (ii) $\hat{SHO}(n,n)_{\epsilon}$ with bracket (5.1) defines a unique non-trivial filtered deformation of $\hat{SHO}(n,n)$.  It is isomorphic to $SHO'(n,n)_{\epsilon}$.

\indent (iii) $\hat{SHO}'(n,n)$ has no simple filtered deformations.

We define a deformed bracket on $SKO(n,n+1;{{{n+2}\over n}})$ as follows (cf.~[Ko]):
\begin{eqnarray}
[f,g]_{\epsilon}&=&([\tau\xi_1\xi_2\cdots\xi_n,fg]+2fg\xi_1\xi_2\cdots\xi_n)\epsilon^{n+2},f,g\in\C[x_1,\cdots,x_n],\nonumber \\
{[f,g]}_{\epsilon}&=&[f,g],\quad{\rm otherwise}.\quad \quad\qquad\;\;\qquad\qquad\qquad\qquad\qquad\qquad{(5.10)}\nonumber 
\end{eqnarray}
Note that $[\tau\xi_1\xi_2\cdots\xi_n,fg]+2fg\xi_1\xi_2\cdots\xi_n\in SKO(n,n+1;{{n+2}\over n})$.

\proclaim{Theorem 5.2} The deformed bracket in (5.10) defines a unique non-trivial filtered deformation of $SKO(n,n+1;{{{n+2}\over n}})$.

\proof  It remains to show that $[\cdot,\cdot]_{\epsilon}$ is a Lie bracket.  Set 
\begin{eqnarray}
\mu_{n+2}(f,g)&=&[\tau\xi_1\xi_2\cdots\xi_n,fg]+2fg\xi_1\xi_2\cdots\xi_n,\quad f,g\in\C[x_1,\cdots,x_n],\nonumber \\
\mu_{n+2}(f,g)&=&0,\quad{\rm otherwise}.\nonumber 
\end{eqnarray}
Note that $\mu_{n+2}(a,\mu_{n+2}(b,c))=0$, for $a,b,c\in SKO(n,n+1;{{{n+2}\over n}})$.  Thus $[\cdot,\cdot]_{\epsilon}$ is a Lie bracket if and only if $\mu_{n+2}$ is a $2$-cocycle of $SKO(n,n+1;{{n+2}\over n})$ with coefficients in its adjoint representa
tion. Hence for $f,g,h\in SKO(n,n+1;{{{n+2}\over n}})$ we only need to check identity (5.4).  It is easy to verify that unless $f$, $g$ and $h$ belong to the $gl_n$-components generated by the highest weight vectors $x_1^k$, $x_1^{k}\xi_n$ and $x_1^k(\tau
+{{\beta n-k}\over{n+k}}\Phi)$ (see the table of $SKO(n,n+1;\beta)$ in Section 4) (5.4) gives the trivial identity $0=0$.  Also it can be verified directly that if two of the $f,g,h$ are not in $\C[x_1,\cdots,x_n]$, then (5.4) again gives the trivial iden
tity.  Hence we are to consider three cases.  Namely

\indent (1) $f,g,h\in\C[x_1,\cdots,x_n]$.

\indent (2) $f,g\in\C[x_1,\cdots,x_n]$ and $h$ is in the $gl_n$-component generated $x_1^k\xi_n$.

\indent (3) $f,g\in\C[x_1,\cdots,x_n]$ and $h$ is in the $gl_n$-component generated by \newline \indent $x_1^k(\tau+{{\beta n-k}\over{n+k}}\Phi)$.

Even though we don't have the Poisson bracket at our disposal, we have the following useful identity:
$$\mu_{n+2}(f,g)=[\tau\xi_1\xi_2\cdots\xi_n,f]g+f[\tau\xi_1\xi_2\cdots\xi_n,g],\quad f,g\in\C[x_1,\cdots,x_n].$$
Note also that $[\tau\xi_1\xi_2\cdots\xi_n,f]$=0, if $f$ lies in the $gl_n$-component generated by $x_1^k\xi_n$ and 
$[\tau\xi_1\xi_2\cdots,\xi_n,g]+2g\xi_1\xi_2\cdots\xi_n=0$, if $g$ lies in the $gl_n$-component generated by $x_1^k(\tau+{{\beta n-k}\over{n+k}}\Phi)$.
Using these identities the computation is similar to the one given in the proof of Proposition 5.1.  {$\qed$\medskip}

\bigskip
\bigskip
\centerline  {\bf REFERENCES }
\bigskip
\frenchspacing
\medskip

\noindent [ALS] Alekseevsky, D.; Leites, D.; Shchepochkina, I.: {\it Examples of simple Lie superalgebras of vector fields}, C.R.~Acad.~Bulg.~Sci.~{\bf 33}, no.~9 (1980) 1187--1190 (in Russian).

\noindent [CK] Cheng, S.-J.; Kac, V.~G.: {\it Structure of some $\Z$-graded Lie Superalgebras of Vector Fields}, preprint.

\noindent [F] Fuchs, D. B.: {\it Cohomology of infinite-dimensional Lie algebras}. New York and London: Consultants Bureau 1986.

\noindent [G] Guillemin, V.~W.: {\it Infinite-dimensional primitive Lie algebras}, J.~Diff.\newline ~Geom.~{\bf 2} (1970) 257--282.

\noindent [K1] Kac, V.~G.: {\it Lie superalgebras}, Adv.~Math.~{\bf 26} (1977) 8--96.

\noindent [K2] Kac, V.~G.: {\it Description of filtered Lie algebras associated with graded Lie algebras of Cartan type}, Math.~USSR Izv.~{\bf 8} (1974) 801--835.

\noindent [K3] Kac, V.~G.: {\it Classification of infinite-dimensional simple linearly compact Lie superalgebras}, Adv.~Math.~{\bf 139} (1998) 1--55.

\noindent [KN] Kobayashi, S.; Nagano, T.: {\it On Filtered Lie Algebras and Geometric Structures, IV}.  J.~of Math.~and Mech.~{\bf 15} (1966) 163--175. 

\noindent [Ko] Kotchetkoff, Yu.: {\it D\'eformation de superalg\'ebras de Buttin et quantification}, C.R.~Acad. Sci.~Paris {\bf 299}, ser.~I, no.~14 (1984) 643--645.

\noindent [L] Leites, D.: {\it Quantization. Supplement 3}. In: F.~Berezin, M.~Shubin:\newline {\it Schr\"odinger equation}, Kluwer, Dodrecht (1991) 483--522.

\noindent [R] Rim, D.~S.: {\it Deformation of transitive Lie algebras}, Ann.~of Math.\newline ~{\bf 83} (1966) 339--357.

\noindent [S] Shchepochkina, I.: {\it The five exceptional simple Lie superalgebras of vector fields}, hep-th/9702121 (1997).

\noindent [Sp] Spencer, D.~C.: {\it Deformation of structures on manifolds defined by transitive continuous pseudogroups}, Ann.~of Math.~{\bf 76} (1962) 306--445.

\noindent [SS] Singer, I.~S.; Sternberg, S.: {\it On the infinite groups of Lie and Cartan I}, J.~Analyse Math.~{\bf 15} (1965) 1--114.

\noindent [W] Weisfeiler, B.~Yu.: {\it Infinite-dimensional filtered Lie algebras and their connection with graded Lie algebras}, Funct.~Anal.~Appl.~{\bf 2} (1968) 88--89.

\end{document}